\documentclass[10pt]{article}
\usepackage{repstyle}
\usepackage{latexsym,epsf,amssymb,subeqn}
\usepackage{graphicx}

\newcommand{\eq}{\begin{equation}}
\newcommand{\eeq}{\end{equation}}
\newcommand{\R}{\mathbb R}
\newcommand{\C}{\mathbb C}
\newcommand{\Real}{{\rm Re}}

\newcommand{\Oh}{{\cal O}}
\newcommand{\sfrac}[2]{\mbox{\large{$#1\over#2$}}}

\newcommand{\sDelta}{{\mbox{\footnotesize{$\Delta$}}}}
\newcommand{\dt}{\sDelta t}
\newcommand{\dx}{\sDelta x}
\newcommand{\dy}{\sDelta y}
\newcommand{\e}{\mathrm e}
\newcommand{\eps}{\varepsilon}
\newcommand{\fhi}{{\varphi}}

\begin{document}

\title{Modified Douglas Splitting Methods for Reaction-Diffusion Equations}

\author{A.\,Arrar\'as%
\footnote{Departamento de Ingenier\'{i}a­ Matem\'{a}tica e Inform\'{a}tica,
Universidad P\'{u}blica de Navarra, Campus de Arrosad{\'i}a,
31006 Pamplona, Spain.
E-mail: andres.arraras@unavarra.es }
, K.J.\,in\,'t\,Hout%
\footnote{
Department of Mathematics and Computer Science, University of Antwerp,
Middelheimlaan 1, B-2020 Antwerp, Belgium.
Email: karel.inthout@uantwerp.be.}
, W.\,Hundsdorfer%
\footnote{CWI, Science Park 123, Amsterdam, The Netherlands.
E-mail: willem.hundsdorfer@cwi.nl}
, L.\,Portero%
\footnote{Departamento de Ingenier\'{i}a­ Matem\'{a}tica e Inform\'{a}tica,
Universidad P\'{u}blica de Navarra, Campus de Arrosad{\'i}a,
31006 Pamplona, Spain.
Email: laura.portero@unavarra.es }
}
\date{}
\maketitle

\begin{abstract} \noindent
We present modifications of the second-order Douglas stabilizing
corrections method, which is a splitting method based on the implicit
trapezoidal rule.
Inclusion of an explicit term in a forward Euler way is straightforward,
but this will lower the order of convergence.  In the modifications
considered here, explicit terms are included in a second-order fashion.
For these modified methods, results on linear stability and convergence
are derived.
Stability holds for important classes of reaction-diffusion equations,
and for such problems the modified Douglas methods are seen to be often 
more efficient than related methods from the literature.
\end{abstract}

\section{Introduction}
\label{Sect:Intro}

In this paper we consider numerical methods for solving
initial value problems for ordinary differential equations
(ODEs) on a time interval $[0,T]$,
\eq
\label{eq:ODE}
u'(t) = F(t,u(t)) \,, \qquad u(0) = u_0 \,,
\eeq
with given $F:\R\times\R^m\rightarrow\R^m$ and $u_0\in\R^m$.
In our applications, this system will be obtained from spatial
discretization of a partial differential equation (PDE) by finite
differences or finite elements (with mass lumping).
Then (\ref{eq:ODE}) is called a semi-discrete system and $m$
will be proportional to the number of grid points in space.
In this formulation, inhomogeneous boundary values for the PDE
are directly incorporated inside $F$ as source terms.

For many problems occurring in practice there will be a natural decomposition
\eq
\label{eq:Deco}
F(t,u) \,=\, F_0(t,u) + F_1(t,u) + \cdots + F_s(t,u)
\eeq
in which the separate component functions $F_j$ are more simple than
the whole $F$. It will be assumed that $F_0$ is a non-stiff or mildly
stiff term that can be treated explicitly in a time stepping method.
The other terms $F_j$ will be treated in an implicit fashion, in such
a way that in each internal stage at most one of these $F_j$ appears
implicitly.

\subsection{Modified Douglas methods}

In the following we denote by $\dt > 0$ the step-size and $u_n \approx
u(t_n)$ stands for the numerical approximation to the exact solution at
time level $t_n = n\dt$, $n\ge0$.

A simple method that employs the splitting (\ref{eq:Deco}) is the
Douglas method, which has been presented in \cite{Hu02,HuVe03} as
\eq
\label{eq:SC1}
\setlength{\arraycolsep}{1mm}
\left\{
\begin{array}{ccl}
v_{0} &=& u_{n-1} + \dt\,F(t_{n-1}, u_{n-1}) \,,
\\[2mm]
v_{j} &=& v_{j-1}
+ \theta\dt \big(F_j(t_{n}, v_{j}) - F_j(t_{n-1}, u_{n-1})\big)
\qquad (j = 1,2,\ldots,s) \,,
\\[2mm]
u_{n} &=& v_{s} \,,
\end{array}
\right.
\eeq
with parameter $\theta \ge {1\over2}$ and with $n=1,2,\ldots$
indicating the time level.
The internal vectors $v_{j}$, which depend of course also on $n$,
are all consistent with the exact solution at time $t_{n}$,
and consequently steady state solutions $F(u)=0$ for autonomous
equations are stationary points of the method.
This is an advantage over many other splitting methods, based on Lie or
Strang splitting, where such consistency is often absent.

Originally, for $F = F_1 + \cdots + F_s$ without explicit terms,
such methods were introduced by J.\,Douglas~Jr and co-workers for
multi-dimensional parabolic problems with dimensional splittings; see
for instance \cite{Do62,DoGu64}, and also \cite[p.\,373]{HuVe03} and
\cite{Ma90} for additional references.
In the terminology of \cite{Ma90}, methods like (\ref{eq:SC1})
are also known as stabilizing correction methods.

If  $\theta = {1\over2}$ and $F_0 = 0$ the method (\ref{eq:SC1}) is
convergent of order two in the classical ODE sense, that is, for fixed
(non-stiff) ODE systems.
However, for non-zero $F_0$ the order will only be one, because the explicit
term is treated in a forward Euler fashion.

To improve this treatment of the explicit terms, we consider the
following modification
\eq
\label{eq:SC1a}
\setlength{\arraycolsep}{1mm}
\left\{
\begin{array}{ccl}
v_* &=& u_{n-1} + \dt\,F(t_{n-1}, u_{n-1}) \,,
\\[2mm]
v_{0} &=& v_* + \sfrac{1}{2}\dt
\big(F_0(t_{n}, v_*) - F_0(t_{n-1}, u_{n-1})\big) \,,
\\[2mm]
v_{j} &=& v_{j-1} + \sfrac{1}{2}\dt
\big(F_j(t_{n}, v_{j}) - F_j(t_{n-1}, u_{n-1})\big)
\qquad (j = 1,2,\ldots,s) \,,
\\[2mm]
u_{n} &=& v_{s} \,.
\end{array}
\right.
\eeq
Now the explicit term is treated as in the explicit trapezoidal rule
(also known as the modified Euler method). As we will see, this method is
indeed of order two in the classical ODE sense, and convergence with
order two will also be valid for interesting classes of PDEs of
reaction-diffusion type.

A closely related variant is obtained if the extra stage with $F_0$
is performed last, leading to
\eq
\label{eq:SC1b}
\setlength{\arraycolsep}{1mm}
\left\{
\begin{array}{ccl}
v_{0} &=& u_{n-1} + \dt\,F(t_{n-1}, u_{n-1}) \,,
\\[2mm]
v_{j} &=& v_{j-1} + \sfrac{1}{2}\dt
\big(F_j(t_{n}, v_{j}) - F_j(t_{n-1}, u_{n-1})\big)
\qquad (j = 1,2,\ldots,s) \,,
\\[2mm]
u_{n} &=& v_{s} + \sfrac{1}{2}\dt
\big(F_0(t_{n}, v_{s}) - F_0(t_{n-1}, u_{n-1})\big) \,.
\end{array}
\right.
\eeq
Again, the explicit term is treated as in the explicit trapezoidal rule,
and the method is of order two in the classical ODE sense. It will be
seen, however, that for semi-discrete systems obtained from PDEs proving
second-order convergence is more difficult than for (\ref{eq:SC1a})
if we have a splitting with $s\ge2$.

These modifications (\ref{eq:SC1a}) and (\ref{eq:SC1b}) were briefly
considered in \cite{Hu02}, but the methods were deemed not to be sufficiently
stable if $F_0$ contains discretized convection terms. In this paper we will
consider applications to reaction-diffusion equations and it will be seen
that for such applications the modified Douglas methods give good
results compared to some other stabilizing correction methods.

For $s=1$, both (\ref{eq:SC1a}) and (\ref{eq:SC1b}) are implicit-explicit
(IMEX) methods where the implicit trapezoidal rule is combined with its
explicit counterpart. These two IMEX methods have been considered and tested
in \cite[pp.\,391--395]{HuVe03}.
Additional splittings of the implicit terms, giving $s \ge 2$, are
important in many practical situations.

\subsection{Outline}

In this paper we will study stability and convergence properties of the
modified Douglas methods (\ref{eq:SC1a}) and (\ref{eq:SC1b}).
The analysis will be performed for linear problems with
$F_j(t,u) = A_j u + g_j(t)$.
Numerical tests will show the relevance for nonlinear problems.

Stability is considered in Section~\ref{Sect:LinStab} in a von Neumann-type
analysis. The results indicate that the modified schemes will be mostly
suited for classes of problems where all implicit terms are known a priori
to have eigenvalues close to the negative real axis, such as
reaction-diffusion equations.

In Section~\ref{Sect:AnalysisErrors} error bounds for the methods are
derived. These error bounds are then discussed in some detail in
Section~\ref{Sect:ModelEqs} for linear
2D model problems $u_t = \Delta u + f$ with dimension splitting.

In Section~\ref{Sect:NumCompar} numerical test results are presented
for nonlinear reaction-diffusion problems.
The modified Douglas methods will be compared to methods from the
literature, in particular from \cite{CrSn88,HoWe09,Hu02}.
Apart from dimension splitting, which requires Cartesian grids,
we will also consider domain decomposition splitting with smooth
overlappings \cite{MPRW98}, which can be used on arbitrary grids in
combination with spatial discretization by finite volume or finite
element methods.

\section{Linear Stability}
\label{Sect:LinStab}

Stability will first be studied for the simple scalar test equation
\eq
\label{eq:TestEq}
u'(t) = \big(\lambda_0 + \lambda_1 + \cdots + \lambda_s\big) u(t) \,,
\qquad \lambda_j \in \C \,.
\eeq
The relevance for PDEs is discussed afterwards.

\subsection{Stability for the scalar test equation}

Consider equation (\ref{eq:TestEq}), and let $z_j = \dt\,\lambda_j$.
Applied to this test equation, both (\ref{eq:SC1a}) and
(\ref{eq:SC1b}) give the recursion $u_n = R\,u_{n-1}$ where
$R = r(z_0,\ldots,z_s)$ is given by the stability function
\eq
\label{eq:StabFunc}
r(z_0,z_1,\ldots,z_s) \,=\, 1 \,+\,
\frac{\big(1 + \frac{1}{2} z_0\big) \, \sum_{j=0}^s z_j}
{\prod_{j=1}^s \big(1 - {1\over2} z_j\big)} \,.
\eeq

In this section we will consider the condition
\eq
\label{eq:StabCond}
|r(z_0,z_1,\ldots,z_s)| \le 1 \,,
\eeq
which provides stability with respect to initial perturbations
on arbitrary time intervals $[0,T]$.
In the results, we will use
$$
{\cal W}_\alpha = \{\zeta \in \C : \; \zeta = 0 \;\; {\rm or} \;\;
|\arg(-\zeta)| \le \alpha \}
$$
to denote closed wedges in the left half-plane, with angle
$\alpha \in [0,\frac{1}{2}\pi]$.

For the case that $z_0 = 0$, $s\ge2$, the following results were
derived in \cite{Hu98}:
\eq
\label{eq:ResA}
|r(0,z_1,\ldots,z_s)\,|\,\le\,1 \;\;
\mbox{for all $z_i \in {\cal W}_\alpha$, $1\le i\le s$}
\quad\Longleftrightarrow\quad
\alpha \,\le\, {1\over s-1} \; {\pi\over2} \;,
\eeq
and if $1\le k\le s-1$ then
\eq
\label{eq:ResB}
\left. \begin{array}{r}
|r(0,z_1,\ldots,z_s)|\,\le\,1 \;\;
\mbox{for all $z_1,\ldots,z_{s-k} \in {\cal W}_\alpha$}
\\
\mbox{and $z_{s-k+1},\ldots,z_{s} \le 0$ } \;\;
\end{array} \right\}
\quad\Longleftrightarrow\quad
\alpha \,\le\,  {1\over s-k} \; {\pi\over2} \;.
\eeq
If $z_0=0$ and $s=1$, then we just have the implicit trapezoidal
rule and stability will hold precisely for $z_1 \in {\cal W}_{\pi/2}$.
Further we note that these results for $z_0=0$ are also relevant if
$|\lambda_0| \le L$ with a fixed $L>0$, corresponding to the case where
$F_0$ is genuinely non-stiff. Then $|z_0| \le \dt\,L$, and with either
condition (\ref{eq:ResA}) or (\ref{eq:ResB}) imposed on $z_1,\ldots,z_s$
we then get $|r(z_0,z_1,\ldots,z_s)| \le 1 + \dt\,K$ with a $K>0$, giving
stability on finite time intervals $[0,T]$.

Of course, for $z_0 = 0$ there is no difference between the
original Douglas method and the modifications (\ref{eq:SC1a}) and
(\ref{eq:SC1b}). For these modifications one would hope to get stability
under the restriction $|1 + z_0 + \frac{1}{2}z_0^2| \le 1$, but
we will see that in general the stronger condition $|1 + z_0| \le 1$
is needed.

The following results provide stability under rather strict conditions.
It will be shown afterwards that these strict conditions are necessary,
that is, the results are sharp.

\begin{Thm} \label{Thm:Stab}
The stability condition (\ref{eq:StabCond}) is valid if either
\eq
\label{eq:stab_s1}
s = 1 \,,  \quad |1 + z_0| \le 1 \,, \quad  z_1 \in {\cal W}_{\pi/2} \,,
\eeq
or
\eq
\label{eq:stab_s2}
s \ge 2 \,, \quad  |1 + z_0| \le 1 \,, \quad  z_j \in {\cal W}_0
\quad (j = 1,2,\ldots,s) \,.
\eeq
\end{Thm}

\medskip\noindent
{\bf Proof.}
For $s = 1$ we have
$$
r(z_0, z_1) \,=\, \frac{ \big(1 + z_0 + \frac{1}{2}z_0^2\big)
+ (1 + z_0) \frac{1}{2}z_1 }{ 1 - \frac{1}{2}z_1 } \;.
$$
According to the maximum modulus theorem it is sufficient to consider
$z_0$ on the boundary of the set where $|1 + z_0| \le 1$. Therefore,
put $1 + z_0 = \e^{i\phi}$ with $i=\sqrt{-1}$. 
Then $1 + z_0 + \frac{1}{2}z_0^2 = \frac{1}{2} (1 + \e^{2i\phi})$. 
Further we can take $z_1$ on the boundary of ${\cal W}_{\pi/2}$, that 
is, $z_1 = i t$ for some $t \in \R$.  For such $z_0, z_1$ we have
$$
|r(z_0, z_1)|^2 \,=\,
\frac{ \frac{1}{4}|1 + \e^{2i\phi}|^2 +
\frac{1}{2} \Real\big( i (1 + \e^{-2i\phi}) \e^{i\phi} \big) t +
\frac{1}{4} |\e^{i\phi}|^2 t^2 }
{1 + \frac{1}{4} t^2} \;.
$$
Since the middle term in the numerator vanishes, this is seen to be
bounded by one.

Next, consider $s\ge2$. If $z_s \rightarrow -\infty$, then
$r(z_0,z_1,\ldots, z_{s-1},z_s)$ converges
to
\eq
\label{eq:r*1}
r_*(z_0, z_1, \ldots, z_{s-1}) \,=\,
1 - 2 \, \frac{\big(1 + \frac{1}{2} z_0\big)}
{\prod_{j=1}^{s-1} \big(1 - \frac{1}{2} z_j\big)} \,.
\eeq
Furthermore, with $\theta_s = (1 - \frac{1}{2} z_s)^{-1}$, we have
\eq
\label{eq:r*2}
r(z_0, \ldots, z_{s-1}, z_s) \,=\,
\theta_s  \, r(z_0, \ldots, z_{s-1})  \,+\,
(1-\theta_s) \, r_*(z_0, \ldots, z_{s-1})  \,,
\eeq
and if $z_s\le0$ (i.e. $z_s\in{\cal W}_0$) then $\theta_s\in[0,1]$.
Hence the maximum of $|r|$ for
$z_s \le 0$ is attained in either $z_s=0$ or $z_s = -\infty$:
\eq
\label{eq:r*3}
\max_{z_s\le0} \, |r(z_0, \ldots, z_{s-1}, z_s)| \,=\,
\max\Big( |r(z_0, \ldots, z_{s-1})|, |r_*(z_0, \ldots, z_{s-1})| \Big)\,.
\eeq

If $|1+z_0|\le1$ and $z_1,\ldots,z_{s-1}\le0$, then
$\theta = \prod_{j=1}^{s-1} (1-{1\over2} z_j)^{-1} \in [0,1]$, and
therefore
$$
|r_*(z_0, \ldots, z_{s-1})| \,=\, |(1-\theta) - \theta (1+z_0)|
\,\le\, (1-\theta) + \theta |1+z_0|
\,\le\, 1 \,.
$$
Sufficiency of (\ref{eq:stab_s2}) now follows easily by induction with
respect to $s$.
\hfill $\Box$

\bigskip
With respect to sharpness of the above results: first note that if $s=1$,
we get, as limit value for $z_1\rightarrow-\infty$,
$$
r_*(z_0) = - (1 + z_0) \,.
$$
This shows that even if $s=1$, then the condition $|1 + z_0| \le 1$ is
already necessary for stability with arbitrary $z_1 \in {\cal W}_0$.
For $s=2$ we get
$$
r_*(z_0, z_1) = - \frac{(1+z_0)+\frac{1}{2}z_1}{1-\frac{1}{2}z_1} \;,
$$
from which it is easily seen that in order to have $|r(z_0,z_1,z_2)| \le 1$
for arbitrary values $|1+z_0|\le1$, $z_1 \in {\cal W}_\alpha$ and
$z_2 \in {\cal W}_0$, we need $\alpha = 0$.

\begin{Rem} \label{Rem:IMEX} \rm
If $s=1$, the above results are also relevant for the Douglas
method (\ref{eq:SC1}) with $\theta=\frac{1}{2}$ because the stability
function of that method equals $-r_*(z_0,z_1)$.
The stability result (\ref{eq:stab_s1}) in Theorem~\ref{Thm:Stab}
for the modified Douglas methods (\ref{eq:SC1a}), (\ref{eq:SC1b})
is more favourable than for the method (\ref{eq:SC1}),
because $|r_*(z_0,z_1)| \le 1$ for all $|1+z_0|\le 1$ will require that
$z_1 \in {\cal W}_0$, instead of $z_1 \in {\cal W}_{\pi/2}$.
\hfill $\Diamond$
\end{Rem}

\subsection{Stability for PDEs}

The above results for the scalar problems are directly applicable to
linear systems with $F_j(u) = A_j u$ and scaled Euclidean norms,
provided the matrices $A_j$ are normal and commuting.
This seems a very strict assumption, but it is common in the analysis
of numerical methods for nonlinear PDEs to use linearization, neglect
boundary conditions and freeze the coefficients in a so-called von
Neumann analysis \cite{Hu02}.
The $\lambda_j$ in the scalar test equation then stand for eigenvalues
of the linearized operators $A_j = \partial F_j(t,u)/\partial u$, and
due to the freezing of coefficients and absence of boundary conditions
these matrices $A_j$ can be assumed to be normal and commuting in the
von Neumann analysis.

Further we note that the above stability conditions have been derived
for the case that the $z_j = \dt \lambda_j$ can be independent
from each other. For convection-diffusion-reaction equations
$u_t + \nabla f(u) = \Delta u + g(u)$, with stiff reaction term $g(u)$,
it might seem attractive to treat the convective term explicitly, say
$F_0(u)\approx \nabla f(u)$, $F_1(u) \approx \Delta u$ and
$F_2(u) \approx g(u)$, and after linearization $z_0$ and $z_1$ are
then related as being eigenvalues for the same Fourier modes.
However, also for this case with related $z_0, z_1$, very strict stability
conditions may be needed.

\begin{figure}[htb]
\setlength{\unitlength}{1cm}
\begin{center}
\begin{picture}(13,4)
\includegraphics[width=13cm]{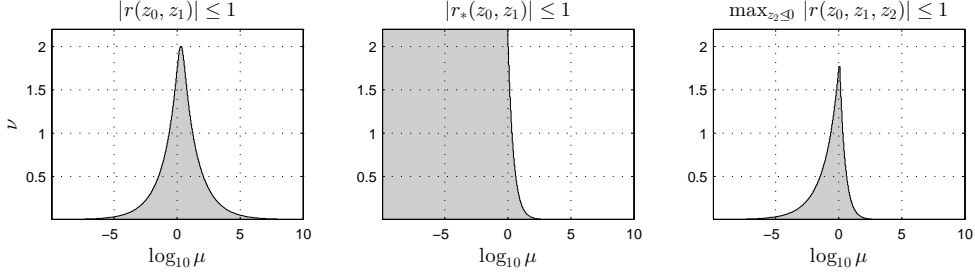}
\end{picture}
\vspace{-3mm}
\caption{ \small  \label{Fig:StReg2}
Domains for stability (gray) in the $(\mu,\nu)$-plane, with cell
P\'eclet number $\mu = a h / d$ and Courant number $\nu = \dt a / h$,
for the modified Douglas schemes with second-order central discretizations.
}
\end{center}
\vspace{-3mm}
\end{figure}

As an example, consider the 1D model problem $u_t + a u_x = d u_{xx} + c u$,
where $c$ is arbitrary negative, representing a stiff reaction term, and
the advection and diffusion terms are discretized with second-order central
differences on a uniform grid with mesh-width $\dx=h$. Then with Courant number
$\nu = \dt\,a/h$ and cell P\'eclet number $\mu = a h / d$, the eigenvalues
are $z_0 = i \nu \sin(2\phi)$, $z_1 = -4 (\nu/\mu) \sin^2(\phi)$ with
$i = \sqrt{-1}$, $\phi\in[0,2\pi]$, and $z_2\le0$ is arbitrary.
It follows by some straightforward (but somewhat tedious) calculations
that $\max_{\phi}|r_*(z_0,z_1)|\le1$ iff $\nu \le 2/\mu$ while the
condition $\max_\phi|r(z_0,z_1)|\le1$  leads to a severe restriction
on the Courant number $\nu$ for small P\'eclet numbers $\mu$
(asymptotically $\nu^3\le4\mu(1+o(\mu))$ for $\mu\rightarrow0$).
The restrictions for stability are displayed in Figure~\ref{Fig:StReg2},
found by direct numerical computations.

If the advection term is discretized by first-order upwind, then we get
eigenvalues $z_0 = - \nu \big( 1 -\cos(2\phi)\big) + i \nu \sin(2\phi)$
which lie on a circle in $\C^-$, with radius $\nu$, touching the imaginary
axis, and therefore $|1 + z_0| \le 1$ if $\nu\le1$.
According to Theorem~\ref{Thm:Stab} this gives stability under
the CFL restriction $\nu\le1$. However, since the goal of the modified
Douglas method is to achieve temporal order two, having only order one
in space is not that interesting.

It seems from the above that advection cannot be taken explicitly with
second-order central spatial discretizations. However, surprisingly, it
can be done if a part of the diffusion term is added to $F_0$.
For example, in 1D, let the discrete advection and diffusion operators
be given, in stencil notation, by
\eq
S_{a} = \sfrac{a}{2 h} \big[ 1\,,\,\,  0\,,\, -1\big] \,, \qquad
S_{d} = \sfrac{d}{h^2} \big[ 1\,,\, -2\,,\,\, 1 \big] \,.
\eeq
Then with $\mu = a h / d$ we can take
\eq
A_0 \,=\, S_{a} + \sfrac{1}{2} \mu S_{d} \,=\,
\sfrac{a}{h} \big[ 1\,,\,  -1\,,\,\, 0\big] \,, \qquad
A_1 \,=\, \big(1 - \sfrac{1}{2} \mu) S_{d} \,,
\eeq
provided the cell P\'eclet number satisfies $\mu \le 2$.
With this $h$-dependent splitting we now get $A_0$ as for first-order upwind,
so stability will hold under the CFL restriction $\nu = \dt\,a/h \le 1$.
In more spatial dimensions a similar splitting can be used.

In conclusion: even though some convection is permitted, the methods
(\ref{eq:SC1a}) and (\ref{eq:SC1b}) with $s\ge2$ seem best suited for
classes of problems where all implicit terms are known a priori to have
eigenvalues close to the real, negative axis.
There are interesting problems of this type originating from systems of
reaction-diffusion equations $u_t = \nabla(D(u)\nabla u) + g(u)$
with diagonal matrix $D$ and with a nonlinear term $g$ which is
either non-stiff or such that the dominant eigenvalues of $g'(u)$
are real negative.
For such problems the methods are expected to be stable.

\section{Analysis of discretization errors}
\label{Sect:AnalysisErrors}

In this section we will analyze the discretization errors for the
modified Douglas methods (\ref{eq:SC1a}), (\ref{eq:SC1b}) in a norm
$\|\cdot\|$ on $\R^m$, with corresponding induced matrix norm on
$\R^{m\times m}$.
The approach will be similar to \cite{Hu02}, where the errors of the
method (\ref{eq:SC1}) have been studied.

We will use the notation $\Oh(\dt^k)$ for a vector or matrix
whose norm is bounded by $C \dt^k$ for small $\dt>0$, with a constant
$C$ that does not depend on the mesh-width $h$ from the spatial
discretization.
The discretization errors will be expressed in terms of the exact
solution $u(t)$ of problem (\ref{eq:ODE}) and
\eq
\label{eq:fhi}
\fhi_j(t) = F_j(t,u(t)) \,.
\eeq
The derivatives of $u(t)$ and $\fhi_j(t)$ appearing in the analysis
are assumed to be $\Oh(1)$ for $t\in[0,T]$, that is, bounded uniformly 
in the mesh-width $h$.

The convergence properties of the methods will be studied for linear
inhomogeneous systems $u'(t) = A u(t) + g(t)$ with
\eq
\label{eq:LinInhom}
F_j(t,v) = A_j v + g_j(t)  \qquad (j = 0,1,\ldots,s)\,.
\eeq
If $A_j$ is a discretized differential operator, with negative powers of
the mesh-width $h$ involved, then inhomogeneous boundary values pertinent
to $A_j$ are supposed to be contained in the source term $g_j(t)$.

In this section it will be assumed that
\eq
\label{eq:assu1}
\big\|\dt A_0 \big\| \le L < 2 \,, \qquad
\big\|(I - \sfrac{1}{2} \dt A_j)^{-1} \big\| \le M \quad (j=1,2,\ldots,s) \,,
\eeq
with a moderate constant $M>0$ (independent of $h$ and $\dt$).
Further we will use the following notation,
$$
Z_j = \dt A_j\,, \quad Z = Z_0 + Z_1 + \cdots + Z_s \,, \quad
P_0 = I + \sfrac{1}{2} Z_0 \,, \quad
Q_j = I - \sfrac{1}{2} Z_j \,.
$$

\subsection{Internal perturbations for method (\ref{eq:SC1a})}

Consider along with (\ref{eq:SC1a}) the scheme with perturbations on each
of the stages
\eq
\label{eq:SC1ap}
\setlength{\arraycolsep}{.5mm}
\left\{
\begin{array}{ccl}
\tilde v_* &=& \tilde u_{n-1} + \dt\,F(t_{n-1}, \tilde u_{n-1})
+ \sigma_0 \,,
\\[2mm]
\tilde v_0 &=& \tilde v_* + \sfrac{1}{2}\dt
\big(F_0(t_{n}, \tilde v_*) - F_0(t_{n-1}, \tilde u_{n-1})\big)
+ \rho_0 \,,
\\[2mm]
\tilde v_{j} &=& \tilde v_{j-1} + \sfrac{1}{2}\dt
\big(F_j(t_{n}, \tilde v_{j}) - F_j(t_{n-1}, \tilde u_{n-1})\big)
+ \rho_{j}
\qquad (j = 1,2,\ldots,s) \,,
\\[2mm]
\tilde u_{n} &=& \tilde v_{s}  \,.
\end{array}
\right.
\eeq
Let $e_n = \tilde u_n - u_n$ and
$w_{j} = \tilde v_{j} - v_{j} - e_{n-1}$,
$w_{*} = \tilde v_{*} - v_{*} - e_{n-1}$.
Subtraction of (\ref{eq:SC1a}) from (\ref{eq:SC1ap}) then gives for
the linear problems (\ref{eq:LinInhom}) the relations
$$
w_* = Z e_{n-1} + \sigma_0 \,, \quad
w_0 = P_0 w_* + \rho_0 \,, \quad
w_{j} = Q_j^{-1}(w_{j-1} + \rho_{j}) \,, \quad
e_n = e_{n-1} + w_{s} \,.
$$
Elimination of the internal quantities $w_*, w_{j}$ leads to the
error recursion
\eq
\label{eq:ErrRec}
e_n \,=\, R\,e_{n-1} + d_n \,,
\eeq
with stability matrix
\eq
\label{eq:StabMat}
R \,=\, I + Q_s^{-1}\cdots Q_2^{-1} \, Q_1^{-1} \, P_0 \, Z
\eeq
and
\eq
\label{eq:LocErr}
d_n \,=\, Q_s^{-1}\cdots Q_2^{-1} Q_1^{-1} \big(P_0 \sigma_0 + \rho_0 \big)
\,+\, \sum_{j=1}^s Q_s^{-1}\cdots Q_j^{-1} \rho_j \,.
\eeq

Under (\ref{eq:assu1}) the norm of $d_n$ can be bounded in terms of the
norms of $\sigma_0$ and the $\rho_j$. To have a favourable error propagation\
for the whole process it will also be assumed that
\eq
\label{eq:assu2}
\|R^n\| \le K  \qquad\mbox{(for all $n\ge1$)} \,,
\eeq
with some moderate constant $K > 0$ (independent of $\dt$ and spatial
mesh-width $h$).
This will hold with $K=1$ in the discrete $L_2$-norm if the matrices
$Z_j$ are normal and commuting, and the eigenvalues $z_j$ satisfy the
linear stability assumptions for the scalar problem (\ref{eq:TestEq}).

\subsection{Local discretization errors for method (\ref{eq:SC1a})}

For the convergence analysis we will use the perturbed method
(\ref{eq:SC1ap}) with $\tilde u_n = u(t_n)$.
Then the $e_n$ in (\ref{eq:ErrRec}) are the global discretization
errors
\eq
\label{eq:GlobErr}
e_n = u(t_n) - u_n \,.
\eeq
Apparently, the term $d_n$ in recursion (\ref{eq:ErrRec}) is then
the local discretization error, introduced in the step from $t_{n-1}$
to $t_n$. To obtain expressions for these local errors, it is convenient
to take $\tilde v_* = \tilde v_{j} =  u(t_n)$. Then the residuals are
\eq
\label{eq:resids}
\begin{array}{c}
\sigma_0 = \sfrac{1}{2}\dt^2 u''(t_{n}) + \Oh(\dt^3) \,,
\\[2mm]
\rho_{j} = \sfrac{1}{2}\dt \big(\fhi_j(t_{n-1}) - \fhi_j(t_{n}) \big)
= -\sfrac{1}{2}\dt^2 \fhi'_j(t_{n}) + \Oh(\dt^3)
\end{array}
\eeq
for $j = 0,1,\ldots,s$.
Since $u''(t) = \sum_{j=0}^s \fhi_j'(t)$,
this gives the following expression for the local errors:
\eq
\label{eq:LocErrSC1a}
\begin{array}{c}
d_n \,=\, \sfrac{1}{2}\dt^2 \, Q_s^{-1}\cdots Q_1^{-1} (P_0 - I)
\big( \fhi'_0(t_{n}) + \fhi'_1(t_{n}) \big)
\\[1mm]
\displaystyle
+\, \sfrac{1}{2}\dt^2 \, Q_s^{-1}\cdots Q_1^{-1}
\sum_{j=2}^s (P_0 - Q_1\cdots Q_{j-1}) \fhi'_j(t_{n})
+ \Oh(\dt^3) \,.
\end{array}
\eeq
It is seen from the derivation that the $\Oh(\dt^3)$ remainder term
can be bounded with a constant determined by the maximum of
$|\fhi_j''(t)|$ for $j=0,1,\ldots,s$ and
$t\in[t_{n-1},t_n]$.

For non-stiff problems, with $Z_j = \Oh(\dt)$, we will have
$d_n = \Oh(\dt^3)$. However, for stiff systems, and in particular
for semi-discrete systems derived from PDEs, we may get a lower
order.
In spite of this, the method will often be second-order convergent
due to damping and cancellation effects, as we will see next.

\subsection{Global discretization errors for method (\ref{eq:SC1a})}

According to the general criterion formulated in \cite{Hu92}, we will have
second-order convergence if the local errors $d_n$ can be decomposed as
\eq
\label{eq:ErrDeco}
d_n \,=\, (R - I) \xi_n + \eta_n \quad \mbox{with} \quad
\xi_n = \Oh(\dt^2), \;
\eta_n = \Oh(\dt^3), \;
\xi_n - \xi_{n-1} = \Oh(\dt^3) \,,
\eeq
and this decomposition was also shown to be necessary in case the $d_n$
are constant, that is, independent of $n$.
We will discuss in detail whether (\ref{eq:ErrDeco}) will be fulfilled.
For simplicity it will be assumed that $Z$ is invertible, although that
is not strictly needed (cf. Remark~\ref{Rem:SingA} below).

To apply criterion (\ref{eq:ErrDeco}), the local errors will be written
in the form
\eq
\label{eq:LocErrZeta}
d_n \,=\, (R - I) \zeta_n + \Oh(\dt^3) \,.
\eeq
These $\zeta_n$ are then natural candidates for $\xi_n$ in (\ref{eq:ErrDeco}).
Furthermore, to express the errors in a compact way, we introduce
\eq
\label{eq:Sj}
S_j \,=\, 2\,(I\,-\,Q_1 Q_2 \cdots Q_j) \qquad (j=1,2,\ldots,s) \,.
\eeq
Then $S_1 = Z_1$, $S_2 = Z_1 + Z_2 - \frac{1}{2} Z_1 Z_2$ and for $2<j\le s$
we have
$$
S_j =\! \sum_{1\le i\le j}Z_i
\,-\!\!\!\! \sum_{1\le i_1<i_2\le j}\!\!\sfrac{1}{2} Z_{i_1} Z_{i_2}
\,+\!\!\!\! \sum_{1\le {i_1}<{i_2}<{i_3}\le j}\!\!\sfrac{1}{4}
Z_{i_1} Z_{i_2} Z_{i_3}
\,+ \cdots +\, \Big(\sfrac{-1}{2}\Big)^{j-1} Z_1 Z_2\cdots Z_j \,.
$$

For the modified Douglas method (\ref{eq:SC1a}) we obtain from
(\ref{eq:LocErrSC1a})
$$
\begin{array}{c}
\displaystyle
d_n \,=\, \sfrac{1}{2}\dt^2 (R-I) Z^{-1} P_0^{-1} \Big[
(P_0 - I) \big( \fhi_0'(t_n) + \fhi_1'(t_n)\big)
\\[2mm]
\displaystyle
+ \, \sum_{j=2}^{s}
(P_0 - Q_1 Q_2 \cdots Q_{j-1}) \fhi'_j(t_n) \Big]
+ \Oh(\dt^3) \,.
\end{array}
$$
Using $P_0 - Q_1\cdots Q_j = \frac{1}{2} (Z_0 + S_j)$, this gives the
expression (\ref{eq:LocErrZeta}) with
\eq
\label{eq:LocErrSC1a'}
\zeta_n \,=\, \sfrac{1}{4}\dt^2 Z^{-1} P_0^{-1} \Big[
Z_0 u''(t_n) + \sum_{j=2}^{s}
S_{j-1} \fhi'_j(t_n) \Big] \,.
\eeq

To find conditions so that $\zeta_n = \Oh(\dt^2)$, we assume for
convenience that $Z_0$ commutes with the other $Z_j$. Then the factor
$P_0^{-1}$ can be moved to the front, and according to (\ref{eq:assu1})
we have $\|P_0^{-1}\| \le (1-\frac{1}{2} L)^{-1}$.
Furthermore, we will have $Z^{-1} S_{j-1} \upsilon(t) = \Oh(1)$ if all
products $Z^{-1} Z_{i_1} \cdots Z_{i_k} \upsilon(t)$ are $\Oh(1)$ for
$1 \!\le\! i_1  \!<\! \cdots  \!<\! i_k  \!<\! j$.
Consequently, criterion (\ref{eq:ErrDeco}) can be applied with
$\xi_n = \zeta_n$ under the following conditions:
\begin{subequations}
\label{eq:Conds}
\eq
\label{eq:Cond1}
A^{-1} A_0 \upsilon(t) \,=\, \Oh(1) \qquad
\mbox{for} \;  \upsilon = u'', u'''
\; \mbox{and} \; t \in [0,T] \,,
\eeq
\eq
\label{eq:Cond2}
\begin{array}{r}
\dt^{k-1} \; A^{-1} A_{i_1} A_{i_2} \cdots A_{i_k} \upsilon(t) \,=\, \Oh(1)
\qquad \mbox{for} \;\; \upsilon = \fhi'_j, \fhi''_j \,,\; t \in [0,T] \;
\mbox{and}
\\
1 \le i_1 < i_2 < \cdots < i_k < j \le s \,.
\end{array}
\eeq
\end{subequations}
In summary, we have obtained the following convergence result.

\begin{Thm} \label{Thm:Conv}
Consider linear problems (\ref{eq:LinInhom}) with
$\fhi_j'(t), \fhi_j''(t) = \Oh(1)$ for $t\in[0,T]$, $j=0,1,\ldots,s$.
Assume (\ref{eq:assu1}) and (\ref{eq:assu2}) hold, $A_0 A_j = A_j A_0$
for $j=1,2,\ldots,s$, and the conditions (\ref{eq:Conds}) are satisfied.
Then method (\ref{eq:SC1a}) will be convergent of order two on the
interval $[0,T]$.
\end{Thm}

If $s=2$, this shows convergence with order two under the condition
(\ref{eq:Cond1}) and $A^{-1} A_1 \upsilon(t) = \Oh(1)$ for
$\upsilon = \fhi'_2, \fhi''_2$, $t \in [0,T]$.
For $s=3$ we get the additional conditions $A^{-1} A_1 \upsilon(t)= \Oh(1)$,
$A^{-1} A_2 \upsilon(t) = \Oh(1)$ and
$\dt \, A^{-1} A_1 A_2 \upsilon(t) = \Oh(1)$
for $\upsilon = \fhi'_3, \fhi''_3$, $t \in [0,T]$.
These results for $s\le3$ are the same as in \cite[Thm.\,2.1, 2.2]{Hu02}
for the case $F_0 = 0$; see also \cite[pp.\,380,\,381]{HuVe03}.
So, in these theoretical results, inclusion
of an explicit term as in the modified Douglas method (\ref{eq:SC1a})
does not lead to conditions for second-order convergence that are more
strict than for $F_0=0$.

The assumption that $A_0$ commutes with the other $A_j$ has been made
only to get more simple conditions. Note that if $A_0 = \Oh(1)$, then
$P_0 = I + \Oh(\dt)$, in which case the conditions (\ref{eq:Conds})
again ensure second-order convergence.

\begin{Exa} \rm
\label{Exa:DimSplit}
The validity of these convergence conditions has been discussed in
\cite{Hu02, HuVe03} for the 3D heat equation in a cube with inhomogeneous
Dirichlet boundary conditions, using dimension splitting for standard
finite differences with uniform mesh-width $h$ in all directions, and
the discrete $L_2$-norm $\|v\|_2 = (\frac{1}{m}\sum_{i=1}^m |v_i|^2)^{1/2}$.
Then $\|A^{-1}\|_2 \le 1$, $\|A^{-1} A_j\|_2 \le 1$. Furthermore, in this
norm, a condition $\dt A^{-1} A_1 A_2 \upsilon = \Oh(1)$ will hold if
$\upsilon$ is the restriction to the grid of a smooth function $w(x)$
which is zero at the boundaries pertinent to either $A_1$ or $A_2$.
If $w$ does not satisfy such a homogeneous boundary condition, we will
still have $\dt A^{-1} A_1 A_2 \upsilon = \Oh(|\log(h)|)$ if $\dt \sim h$
and $\dt A^{-1} A_1 A_2 \upsilon = \Oh(1)$ if $\dt \sim h^{1+\eps}$
for some $\eps>0$.
So in this case convergence with order two is ensured under a mild time
step restriction.

For larger $s$ additional conditions will appear. For example, for $s=4$
the essential new condition will be
$\dt^2 A^{-1} A_1 A_2 A_3 \upsilon(t) = \Oh(1)$
for $\upsilon = \fhi_4', \fhi_4''$, which may lead to additional time
step constraints, depending on the boundary values assumed by the grid
function $\fhi_4$.
\hfill $\Diamond$
\end{Exa}

\begin{Rem} \rm
\label{Rem:SingA}
In the above it has been assumed that $Z = \dt A$ is invertible, but
this is not really needed. The above conditions (\ref{eq:Conds}) are
all of the form $A^{-1} w = \Oh(1)$.
Instead, criterion (\ref{eq:ErrDeco}) can already be applied if
there is a $v \in \R^m$ such that $A v = w$ and $v = \Oh(1)$. This may
hold even if $A$ is singular, as it will be for example for parabolic
problems with Neumann conditions on the whole boundary.
\hfill $\Diamond$
\end{Rem}

\subsection{Discretization errors for method (\ref{eq:SC1b})}

In the same way, we can derive error recursions for method (\ref{eq:SC1b}).
As perturbed scheme we now consider
\eq
\label{eq:SC1bp}
\setlength{\arraycolsep}{.5mm}
\left\{
\begin{array}{ccl}
\tilde v_{0} &=& \tilde u_{n-1} + \dt F(t_{n-1}, \tilde u_{n-1})
+ \sigma_{0} \,,
\\[2mm]
\tilde v_{j} &=& \tilde v_{j-1} + \sfrac{1}{2}\dt
\big(F_j(t_{n}, \tilde v_{j}) - F_j(t_{n-1}, \tilde u_{n-1})\big)
+ \rho_{j}
\qquad (j = 1,2,\ldots,s) \,,
\\[2mm]
\tilde u_{n} &=& \tilde v_{s} + \sfrac{1}{2}\dt
\big(F_0(t_{n}, \tilde v_{s}) - F_0(t_{n-1}, \tilde u_{n-1})\big)
+ \rho_{0} \,.
\end{array}
\right.
\eeq
With $e_n = \tilde u_n - u_n$ and $w_{j} = \tilde v_{j} - v_{j} - e_{n-1}$,
subtraction of (\ref{eq:SC1b}) from (\ref{eq:SC1bp}) gives for
the linear problems (\ref{eq:LinInhom}) the relations
$$
w_{0} = Z e_{n-1} + \sigma_{0}\,, \qquad
w_{j} = Q_j^{-1}(w_{j-1} + \rho_{j}) \,, \qquad
e_n = e_{n-1} + P_0 w_{s} + \rho_{0} \,.
$$
After elimination of the internal quantities $w_{j}$ we obtain the error
recursion (\ref{eq:ErrRec}) with
\eq
\label{eq:StabMatb}
R = I + P_0\,Q_s^{-1}\cdots Q_2^{-1} \, Q_1^{-1} \, Z
\eeq
and
\eq
\label{eq:LocErrb}
d_n = P_0\,Q_s^{-1}\cdots Q_1^{-1}(\sigma_{0}+\rho_{1}) +
P_0\,Q_s^{-1}\cdots Q_2^{-1}\rho_{2} + \cdots +
P_0\,Q_s^{-1}\rho_{s} + \rho_{0} \,.
\eeq

If $\tilde u_n = \tilde v_{j} =  u(t_n)$, the residuals are
the same as in (\ref{eq:resids}).
Using the fact that $u''(t) = \sum_{j=0}^s \fhi_j'(t)$,
this gives the following expression for the local errors:
\eq
\begin{array}{c}
\displaystyle
d_n = \sfrac{1}{2}\dt^2 (P_0\,Q_s^{-1}\cdots Q_1^{-1} - I)
\fhi'_0(t_{n})
\\[2mm]
+\, \sfrac{1}{2}\dt^2 P_0\,Q_s^{-1}\cdots Q_1^{-1}
\sum_{j=2}^s (I - Q_1\ldots Q_{j-1}) \fhi'_j(t_{n})
+ \Oh(\dt^3) \,.
\end{array}
\eeq

By some calculations, it is seen that this local error fits in the
form (\ref{eq:LocErrZeta}) with
\eq
\label{eq:LocErrSC1b'}
\zeta_n \,=\, \sfrac{1}{4}\dt^2 Z^{-1} \Big[
(Z_0 + S_s) P_0^{-1} \fhi_0'(t_n) + \sum_{j=2}^{s}
S_{j-1} \fhi'_j(t_n) \Big] \,.
\eeq
So the main difference between these local errors and the ones for
the method (\ref{eq:SC1a}) is the factor multiplying $\fhi_0'$
in formula (\ref{eq:LocErrSC1b'}). As we will see, this factor will lead
to complications in proving convergence with order two if $s\ge2$.

First consider $s=1$. Then we have the simple expression
$\zeta_n = \frac{1}{4}\dt^2 P_0^{-1} \fhi_0'(t_n)$. Taking $\xi_n = \zeta_n$
thus shows convergence of order two under the sole condition
\eq
\fhi_0'(t), \fhi_0''(t) = \Oh(1) \qquad
\mbox{for} \;t \in [0,T] \,.
\eeq

However, theoretical problems arise for $s=2$. Then we will have a local
error with
$\zeta_n = -\sfrac{1}{8}\dt^2 Z^{-1} Z_1 Z_2 P_0^{-1}\fhi_0'(t_n)
+ \Oh(\dt^2)$,
but in some experiments convergence with order two in the $L_2$-norm was
observed even if this $\zeta_n \neq \Oh(\dt^2)$.
To prove this convergence, a decomposition (\ref{eq:ErrDeco})
with $\xi_n \neq \zeta_n$ will be needed.
As mentioned already, criterion (\ref{eq:ErrDeco}) was shown to be
necessary for constant $d_n$ in \cite[Lemma\,5.2]{Hu92}.
For that proof the following decomposition was used:
\eq
\label{eq:DecoNew}
\xi_n = (\lambda R - I)^{-1} (R - I) \zeta_n \,, \qquad
\eta_n = (\lambda-1) R \, \xi_n + \Oh(\dt^3) \,,
\eeq
with $\lambda = 1/(1\!+\!\dt)$.
In the experiments, this $\xi_n$ was indeed seen to be $\Oh(\dt^2)$.
However, because of the inverse of $\lambda R - I$, this term $\xi_n$ has a
rather complicated structure, and it is not clear how to prove that
it is $\Oh(\dt^2)$ under transparent, realistic conditions.

In conclusion, the variant Douglas method with $s\ge2$ was observed
to be convergent in the $L_2$-norm of order two for many problems, but
simple conditions for such convergence have not been found.
For that, instead of (\ref{eq:DecoNew}) a more tractable decomposition
might be needed.

\begin{Rem} \rm
To prove the necessity of (\ref{eq:ErrDeco}), it was assumed
in \cite{Hu92} that $\|R\|\le1$.
With our stability assumption (\ref{eq:assu2}), where $K\ge1$, this is
achieved with a new norm $\|v\|_* = \sup_{n\ge0}\|R^n v\|$ for $v\in\R^m$.
This new norm is equivalent to the old one, $\|v\| \le \|v\|_* \le K \|v\|$
for all $v \in \R^m$, and we have $\|R\,\|_* \le 1$.
\hfill $\Diamond$
\end{Rem}

\section{Results for 2D model problems with dimension splitting}
\label{Sect:ModelEqs}

\subsection{Convergence in  the discrete $L_2$-norm}

Numerical comparisons of the modified Douglas methods will be presented in
Section~\ref{Sect:NumCompar} for some interesting reaction-diffusion problems.
To illustrate the error build-up, we first show some results for the
Douglas methods applied to a model problem consisting of the inhomogeneous
heat equation
\eq
\label{eq:Model1}
u_t = u_{xx} + u_{yy} + f(x,y,t)
\eeq
with prescribed exact solution
\eq
\label{eq:Model2}
u(x,y,t) = \sin(t)\, \big( (1+2x^2)(1+y^2)-1 \big)
\eeq
on the unit square $(x,y)\in\Omega = [0,1]^2$ with time $t\in[0,1]$ and
Dirichlet boundary conditions.
For this model problem, the source term $f$ as well as the initial and
boundary values are fitted to the exact solution.
Furthermore we use dimension splitting with
standard second-order differences in space and mesh-width $\dx = \dy = h$.
In this example, the source term is put into $F_0$ and $F_1, F_2$ contain
the discretized difference operators in the $x$- and $y$-direction,
respectively, together with the relevant boundary values.

Because the solution is chosen to be a polynomial of degree less than four
in both $x$ and $y$, there will be no spatial errors.  So the errors given
in the following tables are just temporal errors.

In Table~\ref{Tab:TestModLoc} the local errors $\|e_1\|_2 = \|d_1\|_2$
of the first time step are given in the discrete $L_2$-norm for various
step-sizes $\dt$ with the ratio $\dt/h$ held fixed.
For these local errors we see an $\Oh(\dt^2)$ behaviour with
all three schemes. For the original Douglas method (\ref{eq:SC1})
this is in line with it being a first-order method, with an error
proportional to $\dt^2$ after one step. However, for the two modified
methods (\ref{eq:SC1a}) and (\ref{eq:SC1b}) we see a similar behaviour.
This is due to order reduction; on a fixed spatial grid the local errors
are proportional to $\dt^3$, even though very small time steps will be
needed to see that asymptotic behaviour.

\begin{table}[htb]
\vspace{-2mm}
\caption{ \small \label{Tab:TestModLoc}
Local errors $\|e_1\|_2 = \|d_1\|_2$, at time $t_1 = \dt$, for the
heat equation (\ref{eq:Model1}), (\ref{eq:Model2}), with the Douglas
schemes, $\dt/h = 1$.
}
\vspace{-2mm}
\small
\begin{center}
\setlength{\tabcolsep}{3mm}
\begin{tabular}{|c||c|c|c|c|}
\hline \rule[-1mm]{0mm}{5mm}
$1/\dt$ & 50    & 100     & 200     & 400
\\ \hline \hline
\rule{0cm}{4mm}
Meth.\,(\ref{eq:SC1})
 & $1.31\cdot10^{-3}$ & $3.58\cdot10^{-4}$
 & $9.54\cdot10^{-5}$ & $2.49\cdot10^{-5}$
 \\[1mm] \hline
\rule{0cm}{4mm}
Meth.\,(\ref{eq:SC1a})
 & $2.14\cdot10^{-4}$ & $4.91\cdot10^{-5}$
 & $1.10\cdot10^{-5}$ & $2.40\cdot10^{-6}$
 \\[1mm] \hline
\rule{0cm}{4mm}
Meth.\,(\ref{eq:SC1b})
 & $6.70\cdot10^{-4}$ & $1.42\cdot10^{-4}$
 & $3.02\cdot10^{-5}$ & $6.40\cdot10^{-6}$
 \\[1mm] \hline
\end{tabular}
\end{center}
\end{table}

The global errors $\|e_N\|_2$ in the discrete $L_2$-norm at the output
time $t_N=1$ are shown in Table~\ref{Tab:TestModGlob}. Here the Douglas
method (\ref{eq:SC1}) is seen to be convergent with order one,
whereas for the modified methods (\ref{eq:SC1a}) and (\ref{eq:SC1b}) we
get convergence with order two. For method (\ref{eq:SC1a}) this
is in agreement with the result of Theorem~\ref{Thm:Conv}.

\begin{table}[t]
\vspace{-2mm}
\caption{ \small \label{Tab:TestModGlob}
Global errors $\|e_N\|_2$ at final time $t_N = 1$ for the heat equation
(\ref{eq:Model1}), (\ref{eq:Model2}), with the Douglas schemes,
$\dt/h = 1$.
}
\vspace{-2mm}
\small
\begin{center}
\setlength{\tabcolsep}{3mm}
\begin{tabular}{|c||c|c|c|c|}
\hline \rule[-1mm]{0mm}{5mm}
$1/\dt$  & 50    & 100     & 200     & 400
\\ \hline \hline
\rule{0cm}{4mm}
Meth.\,(\ref{eq:SC1})
 & $2.52\cdot10^{-3}$ & $1.22\cdot10^{-3}$
 & $6.04\cdot10^{-4}$ & $3.00\cdot10^{-4}$
 \\[1mm] \hline
\rule{0cm}{4mm}
Meth.\,(\ref{eq:SC1a})
 & $1.21\cdot10^{-4}$ & $3.04\cdot10^{-5}$
 & $7.64\cdot10^{-6}$ & $1.91\cdot10^{-6}$
 \\[1mm] \hline
\rule{0cm}{4mm}
Meth.\,(\ref{eq:SC1b})
 & $6.63\cdot10^{-4}$ & $1.60\cdot10^{-4}$
 & $3.90\cdot10^{-5}$ & $9.60\cdot10^{-6}$
 \\[1mm] \hline
\end{tabular}
\end{center}
\end{table}

For method (\ref{eq:SC1b}) simple convergence criteria are lacking,
but in the table an order two convergence behaviour is observed.
In the present test, the ratio $\dt/h$ is held fixed, and then it can
be shown (cf. Example~\ref{Exa:DimSplit}) that
$\zeta_n = \Oh(|\log(h)|\dt^2)$, which shows convergence with
'practically' order two. If $\dt$ would tend to zero much slower than $h$,
say $\dt\sim\sqrt{h}$, a more complicated error decomposition, such as
(\ref{eq:DecoNew}), would be needed to demonstrate the order two convergence.

Finally we mention that in this test the error constants are smaller
for the modification (\ref{eq:SC1a}) than for the variant (\ref{eq:SC1b}).
It will be seen in Section~\ref{Sect:NumCompar}, with more interesting test 
problems, that this is not always the case.

\subsection{Convergence in the maximum norm}

The global errors for the model problem (\ref{eq:Model1}), (\ref{eq:Model2})
in the maximum norm $\|\cdot\|_\infty$ are given in 
Table~\ref{Tab:TestModGlobInfty}.
As seen from this table, there is now a clear distinction between the
two methods (\ref{eq:SC1a}) and (\ref{eq:SC1b}). Method (\ref{eq:SC1a})
retains its second-order convergence, but the order of convergence for
(\ref{eq:SC1b}) drops to one. In fact, the errors for the latter method
are now even larger than for the original Douglas method (\ref{eq:SC1}).

\begin{table}[t]
\vspace{-2mm}
\caption{ \small \label{Tab:TestModGlobInfty}
Global errors $\|e_N\|_\infty$ at final time $t_N = 1$ for the heat equation
(\ref{eq:Model1}), (\ref{eq:Model2}), with the Douglas schemes,
$\dt/h = 1$.
}
\vspace{-2mm}
\small
\begin{center}
\setlength{\tabcolsep}{3mm}
\begin{tabular}{|c||c|c|c|c|}
\hline \rule[-1mm]{0mm}{5mm}
 $1/\dt$  & 50    & 100     & 200     & 400
\\ \hline \hline
\rule{0cm}{4mm}
Meth.\,(\ref{eq:SC1})
 & $4.37\cdot10^{-3}$ & $2.16\cdot10^{-3}$
 & $1.07\cdot10^{-3}$ & $5.36\cdot10^{-4}$
 \\[1mm] \hline
\rule{0cm}{4mm}
Meth.\,(\ref{eq:SC1a})
 & $3.11\cdot10^{-4}$ & $7.93\cdot10^{-5}$
 & $2.00\cdot10^{-5}$ & $5.04\cdot10^{-6}$
 \\[1mm] \hline
\rule{0cm}{4mm}
Meth.\,(\ref{eq:SC1b})
 & $1.05\cdot10^{-2}$ & $5.04\cdot10^{-3}$
 & $2.46\cdot10^{-3}$ & $1.21\cdot10^{-3}$
 \\[1mm] \hline
\end{tabular}
\end{center}
\end{table}

To understand these results, let us first note that the local errors are
$\Oh(\dt^2)$ for all the methods, just as in Table~\ref{Tab:TestModLoc}.
Furthermore, the methods are stable in the maximum norm.
For this model problem with $A_0=0$ and commuting $A_1$ and $A_2$,
we have $R = R_1 R_2$ with matrices
$R_j = (I-\frac{1}{2}Z_j)^{-1} (I+\frac{1}{2}Z_j)$ for $j=1,2$.
It is known, for the one-dimensional heat equation, that the estimate
$\sup_{n\ge0}\|R_j^n\|_\infty\le \kappa$ is valid with a modest value
$\kappa>0$, cf.\ \cite{BDS02,FaPa02,Ser64}.
Since $R_1$ and $R_2$ commute, it follows that the stability condition
(\ref{eq:assu2}) will hold with constant $K = \kappa^2$.

The local errors for method (\ref{eq:SC1a}) are given by (\ref{eq:LocErrZeta})
with $\zeta_n = \frac{1}{4}\dt^2 A^{-1} A_1 \fhi_2'(t_n)$, and $\fhi_2(t)$
is the restriction to the grid of $u_{yy}(x,y,t)$. In this test we observed
that $\|A^{-1} A_1\|_\infty \sim |\log(h)|$. This leads directly to the
global error bound $\|\e_N\|_\infty = \Oh(|\log(h)|\dt^2)$,
which is 'practically' order two for fixed ratios $\dt/h$.

To show that method (\ref{eq:SC1b}) is only convergent with order one in
the maximum norm, it is convenient to consider the Dirichlet problem for
the heat equation (\ref{eq:Model1}) on the unit square with solution
\eq
\label{eq:Model3}
u(x,y,t) = 1 - \sfrac{1}{2} t x^2 \,.
\eeq
This leads to very simple expressions for the discretization errors,
due to the fact that we now have $u_{tt} = u_{yy} = 0$ and $u_{xx} = -t$
is constant in space. In the following we consider again a uniform grid
with $\dx = \dy = h$, with $k$ points per direction, $h = 1/(k+1)$, giving
a total number of grid points $m = k^2$.

We will use Kronecker products for vectors and matrices. Then the restriction
to the grid of a product function $\phi(x)\psi(y)$ can be written as
$b\otimes a \in \R^m$ with $a,b\in\R^k$, $a_i = \phi(i h)$ and
$b_i = \psi(i h)$. Let $e = (1,1,\ldots,1)^T\in\R^k$.
Since there is no spatial error, $\fhi_1(t)$ is the restriction to the 
grid of $u_{xx} = -t$ and $\fhi_2(t) = 0$, and it readily follows that 
the local errors $d_n$ for method (\ref{eq:SC1b}) are constant in $n$, with
\eq
d_n \,\equiv\, d \,=\, \sfrac{1}{2} \dt^2 \big( Q_2^{-1}Q_1^{-1} - I\big) v
\,, \qquad
v \,=\, e\otimes e \,.
\eeq

Akin to the decomposition (\ref{eq:DecoNew}), let $\lambda = 1/(1+\dt)$ and
consider the vector
\eq
\xi \,=\, (\lambda R - I)^{-1} d \,,
\eeq
where $R = I + Q_2^{-1}Q_1^{-1} Z$ according to formula (\ref{eq:StabMatb})
with $Z_0=0$. It is seen by some calculations that
\eq
\xi \,=\, \sfrac{1}{2} \dt^2 (1+\dt) B^{-1} w
\quad \mbox{with} \quad
B \,=\, -\dt\,Q_1 Q_2 + Z \,,
\;\;
w = (I - Q_1 Q_2) v \,.
\eeq

If $D = \mbox{tridiag}(1, -2, 1)\in\R^{k\times k}$ and
$\mu = \frac{1}{2}\dt/h^2$, then $A_1 = h^{-2} I\otimes D$ and
$A_2 = h^{-2} D\otimes I$ are the finite difference approximations to
$\partial^2/\partial x^2$ and $\partial^2/\partial y^2$, respectively,
and
$$
Q_1 = I\otimes (I-\mu D) \,, \qquad
Q_2 = (I-\mu D)\otimes I \,.
$$
Let $\hat e = D e$. Then $\hat e = -(e_1 + e_k)$, with $e_1, e_k$ the
first and last unit vectors in $\R^k$, and
$$
w \,=\, \mu\, e\otimes \hat e + \mu\,\hat e\otimes e
- \mu^2 \hat e\otimes \hat e \,.
$$
The first component of $w$ equals $w_1 = -2\mu - \mu^2$. Hence
$$
\|w\|_\infty \,\ge\, \mu^2 + 2\mu \,.
$$
In fact, a little inspection shows that equality holds.
Further we have
$$
\|B\|_\infty \,\le \dt\|Q_1 Q_2\|_\infty + \|Z\|_\infty
\,\le\, \dt \|I - \mu\,D\|_\infty^{\,2} + 4 \mu\, \|D\|_\infty
\,=\, \dt (1+4\mu)^2 + 16 \mu\,.
$$
Since $\|B^{-1} w\|_\infty \ge \|w\|_\infty/\|B\|_\infty$ it now follows that
\eq
\|\xi\|_\infty \,\ge\, \sfrac{1}{2} \dt^2
\frac{\mu^2 + 2\mu}{\dt(1+4\mu)^2 + 16 \mu} \,.
\eeq

For the case that $\dt = h$, we have $\mu = 1/(2\dt)$, which gives
\eq
\|\xi\|_\infty \,\ge\,
\frac{\dt}{8 ((2+T)^2 + 8)}
\qquad \mbox{whenever} \quad \dt \le T\,.
\eeq
From the proof of \cite[Lemma\,5.2]{Hu92} it is seen that convergence
of order $p$ in the maximum norm implies $\|\xi\|_\infty = \Oh(\dt^p)$,
and therefore we have $p \le1$. By stability and the local error bound
$\|d_n\|_\infty = \Oh(\dt^2)$ we know the order is at least one, showing
that the order of convergence equals $p=1$ for $\dt=\dx=\dy=h$,
$h\rightarrow0$.

This result with order one convergence for solution (\ref{eq:Model3}) is
in accordance with the errors in Table~\ref{Tab:TestModGlobInfty}
for solution~(\ref{eq:Model2}). For the more simple solution
(\ref{eq:Model3}) the local error of method  (\ref{eq:SC1a}) is
$d_n = \frac{1}{2}\dt^2 Q_2^{-1}(Q_1^{-1}-I)\fhi_2'(t_n) = 0$,
so this method is now exact. This reveals a striking difference
between the two methods (\ref{eq:SC1a}) and (\ref{eq:SC1b}).

\section{Numerical comparisons for reaction-diffusion problems}
\label{Sect:NumCompar}

\subsection{Extended stabilizing correction methods}

For some nonlinear reaction-diffusion problems, we will compare the
modified Douglas methods with two extended methods that require
per step (approximately) twice the computational work of (\ref{eq:SC1a})
and (\ref{eq:SC1b}).

The following method was derived by Hundsdorfer \cite{Hu02}, based on a
Rosenbrock-type method with matrix factorization from Verwer et al.\
\cite{VSBH99}. With internal vectors $v^*_{j}$ and $v_{j}$, the
approximations $u_n$ are computed from
\eq
\label{eq:HuVe}
\setlength{\arraycolsep}{1mm}
\left\{
\begin{array}{ccl}
v^*_{0} &=& u_{n-1} + \dt\,F(t_{n-1}, u_{n-1}) \,,
\\[2mm]
v^*_{j} &=& v^*_{j-1}
+ \theta\dt \big(F_j(t_{n}, v^*_{j}) - F_j(t_{n-1}, u_{n-1})\big)
\qquad (j = 1,2,\ldots,s) \,,
\\[2mm]
v_{0} &=& u_{n-1} + \sfrac{1}{2}\dt
\big(F(t_{n-1}, u_{n-1}) + F(t_n,v^*_s) \big) \,,
\\[2mm]
v_{j} &=& v_{j-1}
+ \theta\dt \big(F_j(t_{n}, v_{j}) - F_j(t_{n}, v^*_{s})\big)
\qquad (j = 1,2,\ldots,s) \,,
\\[2mm]
u_{n} &=& v_{s} \,.
\end{array}
\right.
\eeq
The underlying implicit method is $A$-stable for any $\theta\ge\frac{1}{4}$,
and for $\theta = 1\pm\frac{1}{2}\sqrt{2}$ this implicit method has optimal
damping at infinity ($L$-stability).
Furthermore,  method (\ref{eq:HuVe}) is of order two in the ODE sense for 
any value of $\theta$. 
Some stability and convergence results, pertinent to PDEs, can be found in
\cite{HoMi13,HoWe07,HoWe09,HoWy15,Hu02}. 

A related method, due to in\,'t\,Hout \& Welfert \cite{HoWe09}, is given by
\eq
\label{eq:HoWe}
\setlength{\arraycolsep}{1mm}
\left\{
\begin{array}{ccl}
v^*_{0} &=& u_{n-1} + \dt\,F(t_{n-1}, u_{n-1}) \,,
\\[2mm]
v^*_{j} &=& v^*_{j-1}
+ \theta\dt \big(F_j(t_{n}, v^*_{j}) - F_j(t_{n-1}, u_{n-1})\big)
\qquad (j = 1,2,\ldots,s) \,,
\\[2mm]
v_{0} &=& v^*_{0} +
\sfrac{1}{2}\dt \big(F_0(t_{n}, v^*_{s}) - F_0(t_{n-1}, u_{n-1})\big)
\\[1mm]
& & \;\; + \, \big(\sfrac{1}{2}-\theta\big) \dt
\sum_{j=1}^s \big(F_j(t_n, v^*_s) - F_j(t_{n-1}, u_{n-1}) \big) \,,
\\[2mm]
v_{j} &=& v_{j-1}
+ \theta\dt \big(F_j(t_{n}, v_{j}) - F_j(t_{n-1}, u_{n-1})\big)
\qquad (j = 1,2,\ldots,s) \,,
\\[2mm]
u_{n} &=& v_{s} \,.
\end{array}
\right.
\eeq
Here setting $\theta = \frac{1}{2}$ gives the method of Craig \& Sneyd
\cite{CrSn88}.  For any choice of $\theta$, method (\ref{eq:HoWe}) is 
of order two in the ODE sense, and the underlying implicit method is 
the same as for (\ref{eq:HuVe}).
Results on convergence for PDEs can be found in \cite{HoWy14}. Stability
results for parabolic equations with mixed derivatives were obtained in
\cite{HoMi11,HoMi13,HoWe09}.

In the numerical tests the modified Douglas methods (\ref{eq:SC1a}) and
(\ref{eq:SC1b}) will be compared with the methods (\ref{eq:HuVe}) and
(\ref{eq:HoWe}) with parameter value $\theta = 1-\frac{1}{2}\sqrt{2}$,
and for (\ref{eq:HoWe}) we will consider also $\theta=\frac{1}{2}$,
the Craig-Sneyd method.
In the comparisons the results for the Douglas method (\ref{eq:SC1}) with 
$\theta = \frac{1}{2}$ are also included, even though we know already
that this method is only of order one if $F_0 \neq 0$.
The methods (\ref{eq:HuVe}) and (\ref{eq:HoWe}) were also tested with
parameter value $\theta = \frac{1}{4}$, but this produced errors very 
similar to $\theta = 1-\frac{1}{2}\sqrt{2}$, with lines in the error 
plots practically on top of each other; for clarity those results are 
not included in the figures.

In the plots we will indicate the methods (\ref{eq:HuVe}) and
(\ref{eq:HoWe}) with $\theta = 1-\frac{1}{2}\sqrt{2}$ as HV and HW,
respectively. The Craig-Sneyd method,  
(\ref{eq:HoWe}) with $\theta = \frac{1}{2}$, is indicated as CS.

\subsection{A traveling wave problem with dimension splitting}

As a first test we consider a traveling wave problem in two spatial
dimensions on the unit square $(x,y)\in\Omega = [0,1]^2$. For this
problem we use dimensional splitting  with
$A_1\approx{\partial^2}/{\partial x^2}$,
$A_2\approx{\partial^2}/{\partial y^2}$,
where the approximations are obtained with standard second-order finite
differences on a uniform Cartesian grid with $\dx=\dy=h$.
The nonlinear reaction term is treated explicitly.
An exact solution is known, and time-dependent Dirichlet boundary
conditions are fitted to this exact solution.

The equation is given by
\begin{subequations}
\label{eq:travwave1}
\eq
u_t = \epsilon(u_{xx} + u_{yy}) + \gamma u^2 (1-u) \,,
\eeq
with $\gamma, \epsilon>0$, for which we have the traveling wave solution
\eq
u(x,y,t) \,=\, \frac{1}{1+\exp(\beta(r - r_0))} \;,
\qquad
r = \cos(\alpha) x + \sin(\alpha) y - c\,t \,,
\eeq
\end{subequations}
with constant $\beta = \frac{1}{2}\sqrt{{2\gamma}/{\epsilon}}$ and wave 
velocity $c = \sqrt{{\gamma\epsilon}/{2}}$. The angle $\alpha$ and shift 
parameter $r_0$ are arbitrary; in the test we will take 
$\alpha = \frac{1}{6}\pi$ and $r_0 = 1-c$. 
Further we take $\gamma = 50$, and for $\epsilon$ two values
are considered: $\epsilon = 1$ and $\epsilon = \frac{1}{50}$.
For the larger value of $\epsilon$ the solution is very smooth, and the
time-dependent boundary conditions then influence the errors of the schemes.
For $\epsilon = \frac{1}{50}$ the largest errors are found, for all schemes,
in the interior of the domain, showing that the boundary conditions are then
no longer a dominating error source.

A similar problem was used in \cite{Hu02} for tests with the Douglas method
(\ref{eq:SC1}) where the nonlinear reaction term was treated implicitly.
Since this reaction term is not very stiff, it is more natural to treat
it in an explicit fashion, as will be done here. Method (\ref{eq:SC1})
will then be first-order convergent only.

The errors at the output time $T=1$ in the discrete $L_2$-norm and
maximum norm are presented in the plots as function of the mesh-width $h$.
The time step is taken as $\dt=h$ for the Douglas methods (\ref{eq:SC1}),
(\ref{eq:SC1a}), (\ref{eq:SC1b}), and $\dt=2h$ for the extended methods
(\ref{eq:HuVe}) and (\ref{eq:HoWe}). The computational work to cover the
time interval $[0,T]$ is then very similar for all the schemes.
Further it should be noted that the errors in the plots are the total
errors in space and time, but it was verified that the temporal errors
are the dominating ones.

\begin{figure}[t]
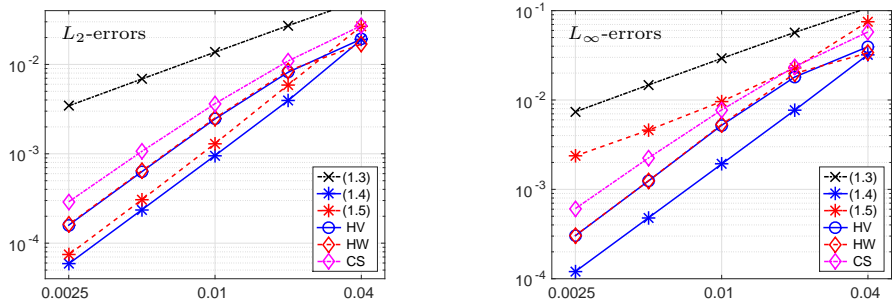

\setlength{\unitlength}{1cm}
\begin{center}
\begin{picture}(5,4.0)
\includegraphics[width=5cm]{FigErrTW1a1.eps}
\put(-4.3,3.5){\scriptsize $L_2$-errors}
\end{picture}
\hspace{1.5cm}
\begin{picture}(5,4.0)
\includegraphics[width=5cm]{FigErrTW1a2.eps}
\put(-4.3,3.5){\scriptsize $L_\infty$-errors}
\end{picture}
\vspace{-3mm}
\caption{ \small  \label{Fig:Resu1a}
$L_2$-errors [left] and $L_\infty$-errors [right] versus mesh-width 
$h = \frac{1}{25},\frac{1}{50},\ldots,\frac{1}{400}$,
for problem (\ref{eq:travwave1}) with $\epsilon = 1$, output time $T = 1$.
Time step $\dt=h$ for the Douglas methods (\ref{eq:SC1}), (\ref{eq:SC1a}), 
(\ref{eq:SC1b}), and $\dt=2h$ for the extended methods
HV [i.e.\ (\ref{eq:HuVe}) with $\theta=1-\frac{1}{2}\sqrt{2}$],
HW [i.e.\ (\ref{eq:HoWe}) with $\theta=1-\frac{1}{2}\sqrt{2}$] and
CS [i.e.\ (\ref{eq:HoWe}) with $\theta=\frac{1}{2}$].
}
\end{center}
\vspace{-3mm}
\end{figure}

The errors for the smooth case, with $\epsilon = 1$, are given in
Figure~\ref{Fig:Resu1a}. As for the linear problem of the previous
section, method (\ref{eq:SC1b}) suffers from order reduction in the
maximum norm, due to the boundary conditions. In the $L_2$-norm the
order of convergence is two with this method, showing that the largest
errors are confined to small regions near the boundaries.
The errors for the methods (\ref{eq:HuVe}) and (\ref{eq:HoWe}) with
$\theta = 1-\frac{1}{2}\sqrt{2}$ are very similar; the errors for
(\ref{eq:HoWe}) are slightly smaller but this is not really visible.
For (\ref{eq:HoWe}) with $\theta=\frac{1}{2}$ (the Craig-Sneyd method) the
errors are a little larger.

\begin{figure}[b!]
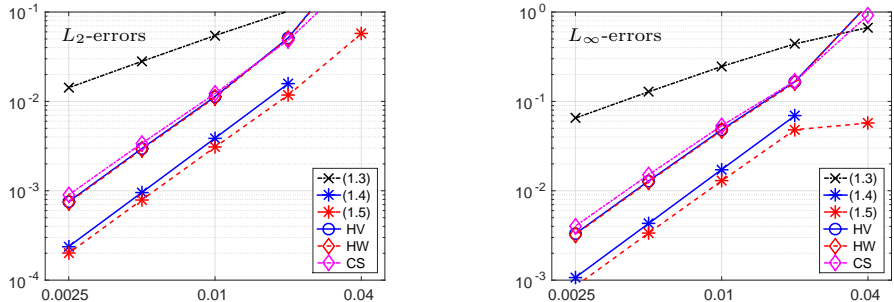

\setlength{\unitlength}{1cm}
\begin{center}
\begin{picture}(5,4.5)
\includegraphics[width=5cm]{FigErrTW1b1.eps}
\put(-4.3,3.5){\scriptsize $L_2$-errors}
\end{picture}
\hspace{1.5cm}
\begin{picture}(5,4.0)
\includegraphics[width=5cm]{FigErrTW1b2.eps}
\put(-4.3,3.5){\scriptsize $L_\infty$-errors}
\end{picture}
\vspace{-3mm}
\caption{ \small  \label{Fig:Resu1b}
$L_2$-errors [left] and $L_\infty$-errors [right] versus mesh-width 
$h = \frac{1}{25},\frac{1}{50},\ldots,\frac{1}{400}$,
for problem (\ref{eq:travwave1}) with $\epsilon = \frac{1}{50}$, output
time $T = 1$.
Time step $\dt=h$ for the Douglas methods, and $\dt=2h$ for the
extended methods (legend names as in Figure~\ref{Fig:Resu1a}).
}
\end{center}
\vspace{-3mm}
\end{figure}

The results for $\epsilon = \frac{1}{50}$ are shown in Figure~\ref{Fig:Resu1b}.
Here method (\ref{eq:SC1b}) is convergent with order two in the
maximum norm. In fact, the errors for this method are now slightly
smaller than for (\ref{eq:SC1a}).
Both modifications of the Douglas method give smaller errors than
the extended stabilizing correction methods (\ref{eq:HuVe}) and
(\ref{eq:HoWe}). These extended methods produced very similar results
in this test.

For $h = \frac{1}{25}$ method (\ref{eq:SC1a}) is unstable, and the same was
observed for the methods (\ref{eq:HuVe}) and (\ref{eq:HoWe}) with
$\theta = \frac{1}{4}$.
Surprisingly, method (\ref{eq:SC1b}) is more stable here.
It should be noted, however, that for $h = \frac{1}{30}$ all
the methods in the test are stable whereas for $h = \frac{1}{20}$
all the methods become unstable, so the differences between the
methods is not very large in this respect.

\subsection{A pattern formation problem on a hexagon with domain
decomposition splitting}

Dimensional splitting requires Cartesian grids. This can be a serious
drawback in applications. Non-Cartesian grids, in combination with
finite elements (with mass lumping) or finite volume discretizations, 
can be used with splittings based on domain decomposition, as
introduced by Mathew, Polyakov, Russo \& Wang \cite{MPRW98}.

\begin{figure}[b!]
\setlength{\unitlength}{1cm}
\begin{center}
\begin{picture}(5,5)
\put(-2.25,0.0){\includegraphics[width=9cm]{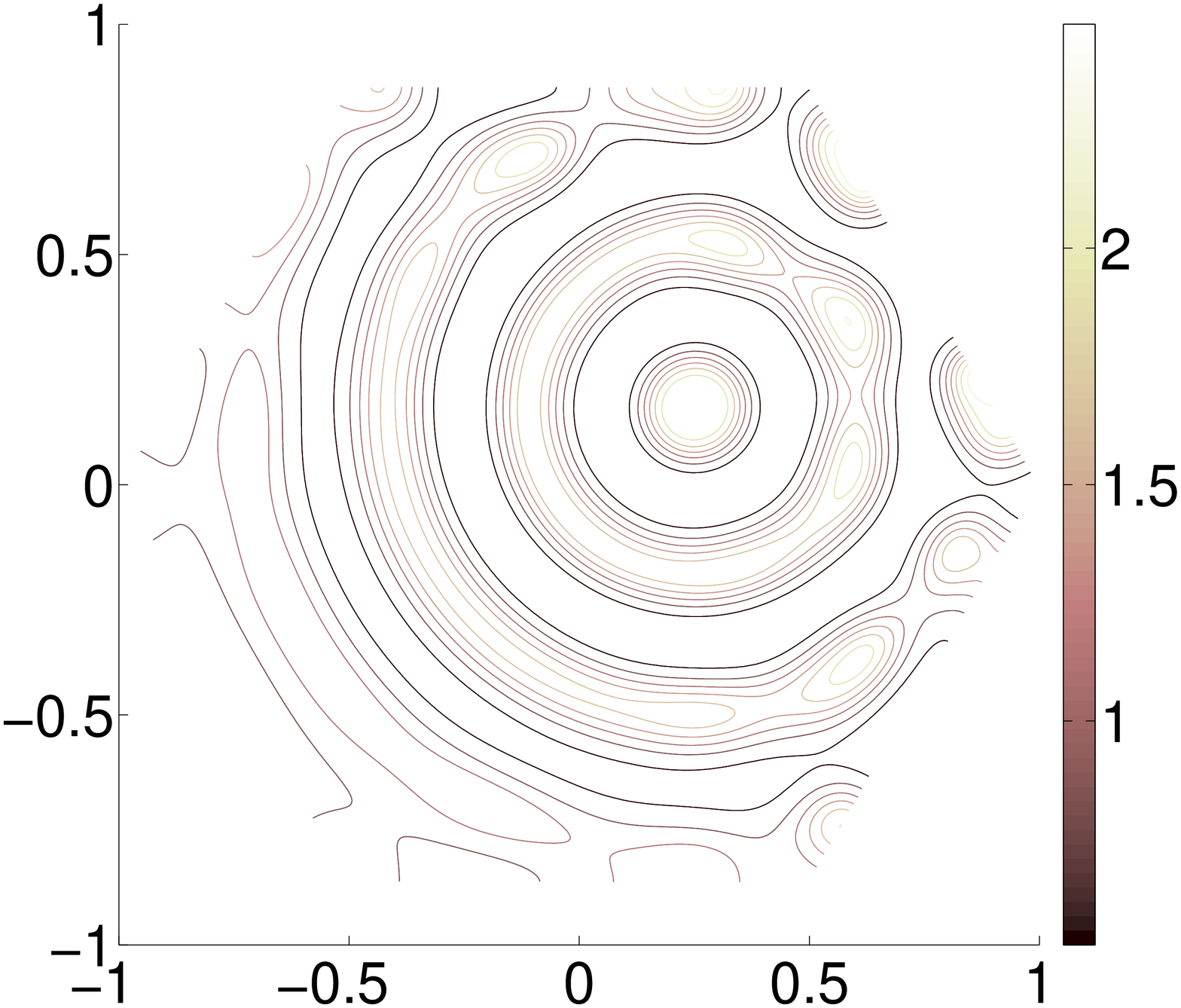}}
\put(0.25,4.4){\scriptsize \fbox{$t=0.5$}}
\end{picture}
\hspace{10mm}
\begin{picture}(5,5)
\put(-2.25,0.0){\includegraphics[width=9cm]{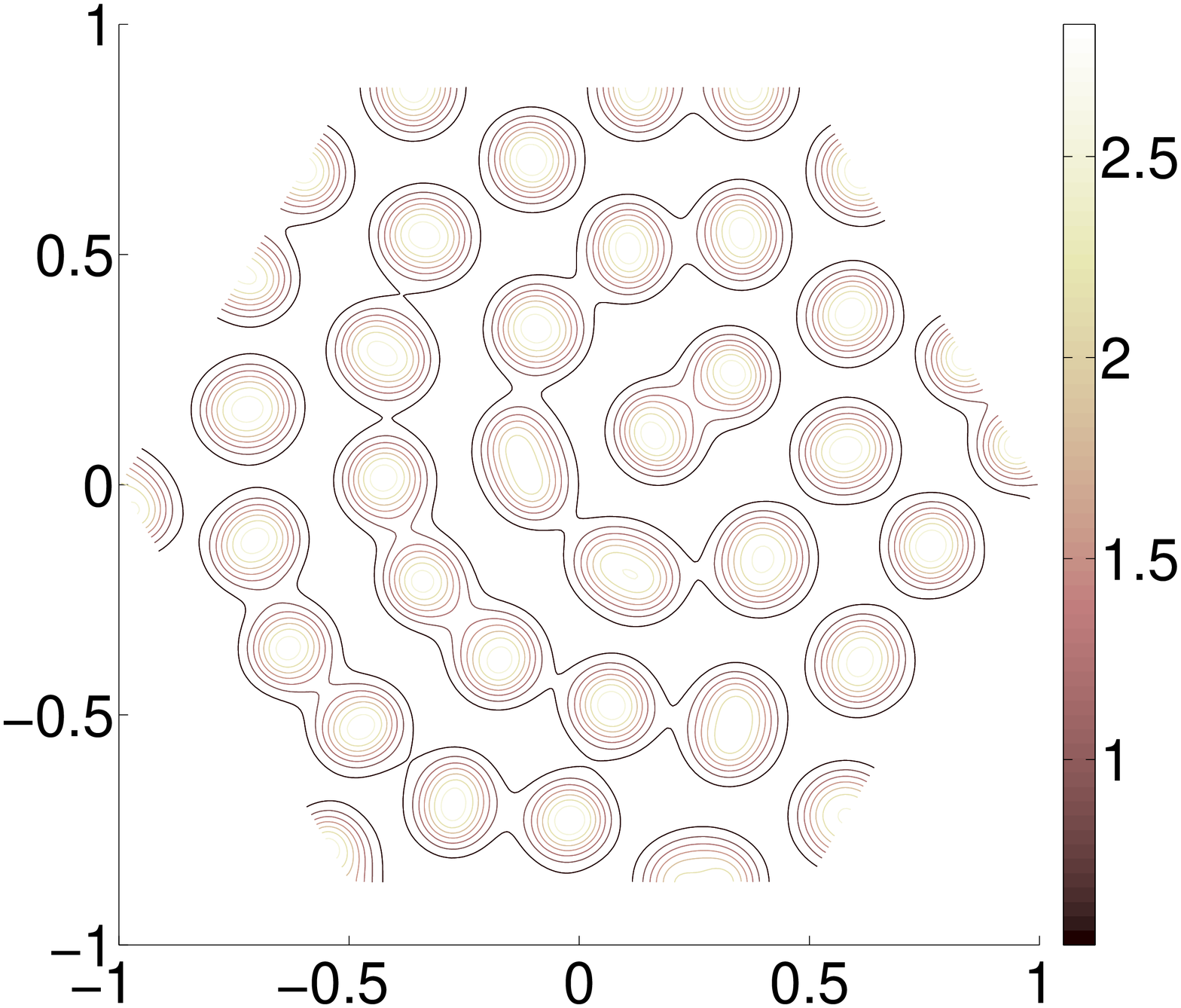}}
\put(0.25,4.4){\scriptsize \fbox{$t=1$}}
\end{picture}
\vspace{-3mm}
\caption{ \small  \label{Fig:SchnackSol}
Contour plots of the $u$-component at times $t=0.5$ [left] and $t=1$ [right].
}
\end{center}
\vspace{-3mm}
\end{figure}

To illustrate this we consider the following system of reaction-diffusion 
equations, describing the interaction of two chemical species,
\eq
\setlength{\arraycolsep}{1mm}
\label{eq:Schnak}
\begin{array}{ccl}
u_t &=& D_1(u_{xx}+u_{yy}) \,+\, \kappa(a-u+u^2v) \,,
\\[2mm]
v_t &=& D_2(v_{xx}+v_{yy}) \,+\, \kappa(b-u^2v) \,,
\end{array}
\eeq
on a hexagonal spatial domain $\Omega$. The vertices of this hexagon are 
located at the points $(\pm1,0)$, $(\frac{1}{2}\pm1,\frac{1}{2}\sqrt{3})$ and
$(\frac{1}{2}\pm1,-\frac{1}{2}\sqrt{3})$. The initial condition is
$$
u(x,y,0) = a + b + 10^{-3}
e^{-100\left(\left(x-\frac{1}{4}\right)^2+\left(y-\frac{1}{6}\right)^2\right)},
\qquad
v(x,y,0) = {b}/{(a+b)^2},
$$
and at the boundaries homogeneous Neumann conditions are imposed. 
The parameter values are $D_1=0.05$, $D_2=1$, $\kappa=100$, $a=0.1305$ 
and $b=0.7695$. The initial condition consists of a small Gaussian 
perturbation added to the chemical steady state $u \equiv a+b = 0.90$, 
$v \equiv b/(a+b)^2 = 0.95$. Due to the reaction and diffusion processes, 
this small perturbation is amplified and spread, leading to the formation 
of patterns with spots. 
Figure~\ref{Fig:SchnackSol} shows snapshots of the time evolution of the 
$u$-component at times $t=0.5$ and $t=1$. The values for $u$ vary between 
$0.2$ and $2.8$, whereas those for $v$ range from $0.4$ to $1.2$, with 
smaller values in the spots where $u$ is maximal.

This model is due to Schnakenberg \cite{S79}. It was considered on a 
rectangular domain in \cite{HuVe03} as a test for IMEX schemes with $s = 1$, 
treating the 2D diffusion terms implicitly and the nonlinear reaction terms 
explicitly. It should be noted that the IMEX method~(4.12)
of \cite[p.\,391]{HuVe03}, which gave favourable results in that test, 
corresponds to the modified Douglas method (\ref{eq:SC1b}) in case $s=1$.

\begin{figure}[b!]
\setlength{\unitlength}{2cm}
\begin{center}
\begin{picture}(3,2.4)
\put(-0.8,0.0){\includegraphics[width=5.5cm]{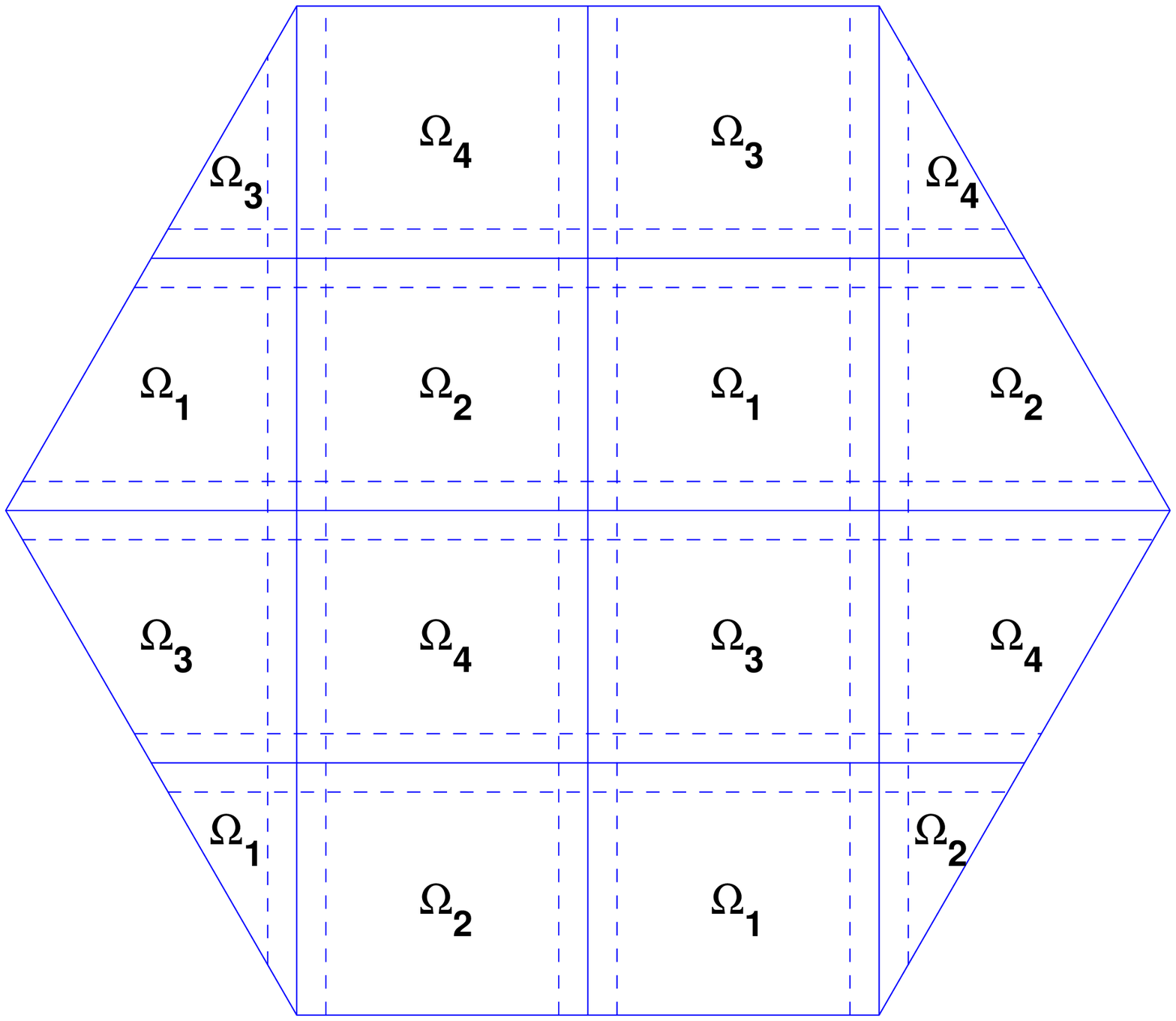}}
\put(1.60,0.1){\scriptsize $s=4$}
\end{picture}
\hspace*{-3cm}
\begin{picture}(3,2.4)
\put(0.8,0.0){\includegraphics[width=5.5cm]{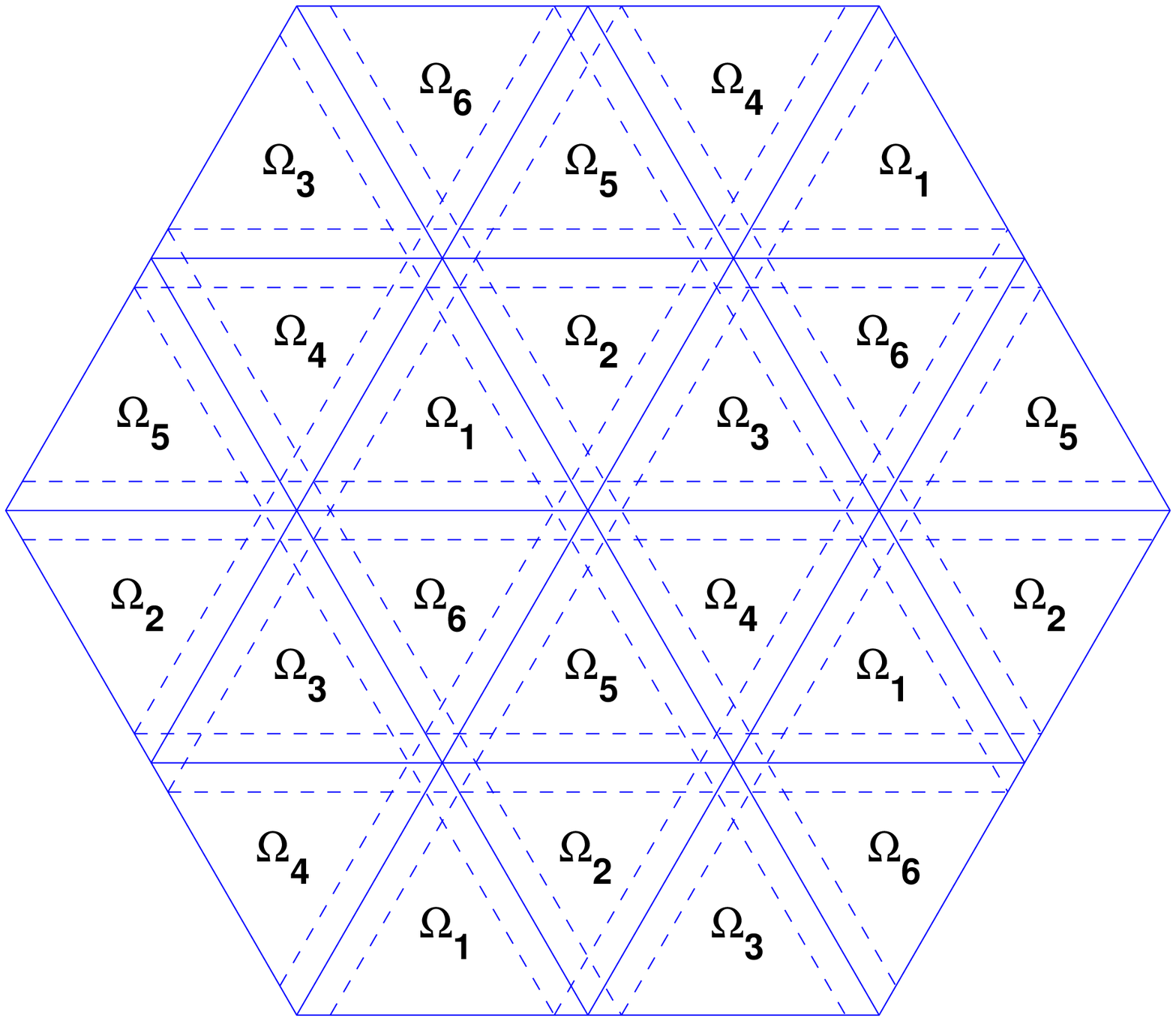}}
\put(3.30,0.1){\scriptsize $s=6$}
\end{picture}
\vspace{-3mm}
\caption{ \small  \label{Fig:HexagDecomp}
Decomposition of the hexagonal domain $\Omega$ into squares, $s=4$ [left] 
and triangles, $s=6$ [right]. The dashed lines indicate the overlap of the
regions.
}
\end{center}
\vspace{-3mm}
\end{figure}

Here, we consider the model to test the behaviour of the stabilizing 
correction methods with $s>2$, using domain decomposition splitting of 
the diffusion terms, discretized by finite elements on a triangular grid.
We will primarily consider domain decomposition splitting with subdomains 
formed by squares restricted to $\Omega$, giving $s=4$, as in 
Figure~\ref{Fig:HexagDecomp}.
In that figure also the subdomains formed by triangles with $s=6$ are shown; 
comments for this are given at the end of this section.
The domain decomposition splitting technique is briefly outlined here
for the squares. A more general and detailed description of this technique 
can be found in \cite{ArPo15,MPRW98}; for further applications we refer to
\cite{AGPR15,APY14}. 

Consider a two-dimensional problem posed on a region $\Omega$ which is 
embedded in a square $[a,b]^2$. Suppose $K$ is a positive integer and
$z_k = a + (b-a)(k-\frac{1}{2})/K$ for $k=1,2,\ldots,K$.
Let $\psi(z):[a,b]\rightarrow[0,1]$ be a smooth function 
(e.g.\ a spline) such that
$$
\psi(z) = \left\{
\begin{array}{ll}
1 & \mbox{if}\;\; |z - z_k|< r \;\;\mbox{with $k$ odd} \,,
\\[1mm]
0 & \mbox{if}\;\; |z - z_k|< r \;\;\mbox{with $k$ even} \,,
\end{array}
\right.
$$
where $r < 1/(2K)$. We now define $\psi_1(x,y) = \psi(x)\psi(y)$,
$\psi_2(x,y) = (1-\psi(x))\psi(y)$, $\psi_3(x,y) = \psi(x)(1-\psi(y))$
and $\psi_4(x,y) = (1-\psi(x))(1-\psi(y))$.
Then, the two-dimensional Laplace operator $\Delta=\nabla^2$ can be split 
in a natural way into four parts, $\Delta=\sum_{j=1}^4\Delta_j$ with
$$
\Delta_j u = \nabla\big(\psi_j(x,y)\nabla u\big)
\quad \mathrm{for}\; 1\leq j\leq 4.
$$
With $\Omega_j$ being the support of $\psi_j$ in $\Omega$, it is seen that 
$\Delta_j u$ will be zero outside $\Omega_j$. Moreover, $\Omega_j$
consists of a number of disconnected component sets. Consequently, if
$A_j$ is the discretization of $\Delta_j$, then a linear system with 
matrix $I - \theta\dt A_j$ will break up into a number of parallel linear 
systems, each corresponding to one of the component sets of $\Omega_j$.

Since the subdomains in which the spatial domain is decomposed
do not have to be aligned with the PDE domain nor with the spatial
discretization grid, this technique is much more flexible in this sense
than the (more classical) alternating direction approach.

Following the above procedure, we consider a decomposition of our hexagonal 
domain $\Omega$ into four overlapping subdomains $\{\Omega_j\}_{j=1}^{4}$, 
each of which consists of four disjoint components, as shown in 
Figure~\ref{Fig:HexagDecomp} (left).

The equations (\ref{eq:Schnak}) are discretized on a triangular grid 
composed of $98304$ equilateral triangles with length $h = \frac{1}{256}$
for the sides, using standard 
piecewise linear finite elements with mass lumping.
On this fixed spatial grid we take time steps 
$\dt = \frac{1}{200}, \frac{1}{400}, \ldots, \frac{1}{6400}$ 
for the Douglas methods (\ref{eq:SC1}), (\ref{eq:SC1a}) and (\ref{eq:SC1b}).
For the extended methods (\ref{eq:HuVe}) and (\ref{eq:HoWe}), the time steps
are taken twice as large, again ensuring that the work to cover a time 
interval $[0,T]$ is the same for all methods.
For these step-sizes the temporal errors for the $u$-component are plotted
in the discrete $L_2$-norm. These errors have been obtained by comparison 
with a time-accurate reference solution on the same grid. 

\begin{figure}[t]
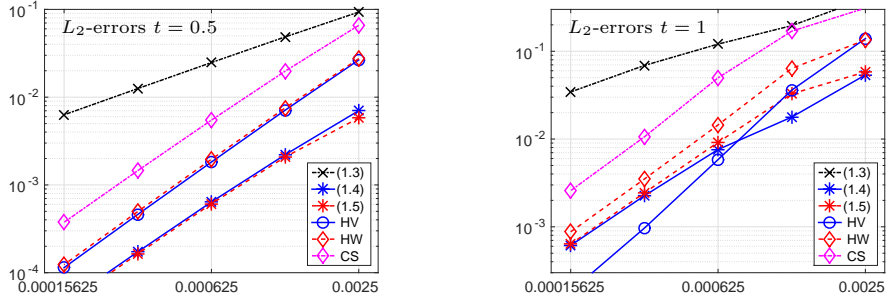

\setlength{\unitlength}{1cm}
\begin{center}
\begin{picture}(5,4.5)
\includegraphics[width=5cm]{FigErrS1.eps}
\put(-4.3,3.5){\scriptsize $L_2$-errors $t=0.5$}
\end{picture}
\hspace{1.5cm}
\begin{picture}(5,4.0)
\includegraphics[width=5cm]{FigErrS2.eps}
\put(-4.3,3.5){\scriptsize $L_2$-errors $t=1$}
\end{picture}
\vspace{-3mm}
\caption{ \small  \label{Fig:SchnackL2Errs}
$L_2$-errors versus $\tau$ for problem (\ref{eq:Schnak}) at
$t=0.5$ [left] and $t=1$ [right] with fixed spatial grid and
$\tau = \frac{1}{400},\frac{1}{800},\ldots,\frac{1}{6400}$.
Step-sizes $\dt = \tau$ for the Douglas methods, and $\dt = 2 \tau$
for the extended methods (legend names as in Figure~\ref{Fig:Resu1a}).
}
\end{center}
\vspace{-3mm}
\end{figure}

The left panel of Figure~\ref{Fig:SchnackL2Errs} shows the errors at 
time $t=0.5$. The Douglas method is seen to be first-order convergent 
(since $F_0\neq 0$), while the other methods all converge with order two.
Among them, the Craig-Sneyd method gives the largest errors, followed by 
the extended methods (\ref{eq:HuVe}) and (\ref{eq:HoWe}), with 
$\theta = 1 - \frac{1}{2}\sqrt{2}$, which perform very similar. 
Finally, the modified Douglas methods (\ref{eq:SC1a}) and (\ref{eq:SC1b}) 
provide the most accurate solutions.

The right panel of Figure~\ref{Fig:SchnackL2Errs} displays the errors 
at time $t=1$. Here, similar comments apply for all the methods, except 
for method (\ref{eq:HuVe}) which is now the most accurate for small 
step-sizes. The reason why this method performs better for this larger 
output time is not clear. It seems that the errors keep increasing in 
time for all the schemes, but this increase is much larger for 
method (\ref{eq:HoWe}) and the modified Douglas methods than for method
(\ref{eq:HuVe}). 
Some additional tests showed that this behaviour is irrespective of the 
value of the parameter $\theta$ in the extended methods.

As seen in Figure~\ref{Fig:HexagDecomp}, the decomposition of the
hexagon with squares gives component sets of different size.
To balance the size of the disjoint components, we could also consider 
products of functions in the directions given by
$y=0$, $y=\sqrt{3}x$ and $y=-\sqrt{3}x$. This would lead to a domain
decomposition into $6$ subdomains, each consisting of disjoint triangular
component sets, as shown in the right panel of
Figure~\ref{Fig:HexagDecomp}, giving a value $s=6$ in the splitting formula
(\ref{eq:Deco}). For this decomposition the errors of the stabilizing 
correction methods were very similar to the decomposition with squares, 
again with a somewhat different behaviour at the different output times, 
so these results are not reproduced here.

\section{Concluding remarks}

As we have seen in the tests for reaction-diffusion problems, the modified
Douglas methods (\ref{eq:SC1a}) and (\ref{eq:SC1b}) can be considerably more
efficient than the extended stabilizing correction methods (\ref{eq:HuVe}) 
and (\ref{eq:HoWe}).

For more general applications, it should be noted that the modified
Douglas methods may not be as stable as the extended methods with
parameter $\theta>\frac{1}{4}$. For example, for parabolic problems with
cross derivatives that are treated explicitly, instabilities have been
observed with the modified Douglas methods while the extended methods still
give good results with properly chosen parameter values $\theta$, see e.g.\ 
\cite{HoMi11}--\cite{HoWy14}.

The variant (\ref{eq:SC1b}) is more likely to suffer from order reduction
than the modification (\ref{eq:SC1a}).
Such an order reduction can be avoided by using boundary corrections, see
e.g.\ \cite{PJB04}. However, since boundary corrections have to be derived
anew for each problem, a method that does not need this has a distinct
advantage over methods that do need such corrections.

\bigskip\noindent
{\bf Acknowledgment}. The work of A. Arrar\'{a}s and L. Portero 
was partially supported by MINECO grant MTM2014-52859.
These authors also gratefully acknowledge the hospitality of the 
Centrum Wiskunde \& Informatica, where this research was partly 
carried out.

\end{document}